\documentclass{article}

\usepackage[T1]{fontenc}
\usepackage[utf8]{inputenc}
\usepackage[english]{babel}

\usepackage{amsmath}
\usepackage{amssymb}
\usepackage{bm}
\usepackage{bbm}
\usepackage{spalign}
\usepackage{enumerate}
\usepackage{fix-cm}
\usepackage{paralist}
\usepackage{enumitem}
\usepackage{multicol}

\usepackage{graphicx} 
\usepackage{xcolor}
\usepackage{tikz}
\usetikzlibrary{positioning, arrows.meta, shapes.geometric}
\usepackage{pgfplots}
\pgfplotsset{compat=1.18}
\usepgfplotslibrary{groupplots}
\usepackage{caption}
\usepackage{float}
\usepackage{subfig}
\usepackage{wrapfig}
\usepackage{longtable}
\usepackage{afterpage}

\usepackage{charter}
\newcommand{\package}{\texttt{Boscia.jl}\xspace}

\usepackage{url} % [obeyspaces]

\graphicspath{{images/}}

\usepackage{fancyvrb}

\IfFileExists{jlcode.sty}{
  \usepackage[theme=default]{jlcode}
}{
  \usepackage{listings}
  \providecommand{\jlbasicfont}{\ttfamily\small\selectfont}
  \lstdefinestyle{jlcodestyle}{
    basicstyle=\jlbasicfont,
    breaklines=true,
    columns=fullflexible,
    keepspaces=true,
    showstringspaces=false
  }
  \NewDocumentEnvironment{jllisting}{o g}{
    \IfNoValueTF{#1}{\lstset{style=jlcodestyle}}{\lstset{style=jlcodestyle,#1}}
    \lstset{basicstyle=\jlbasicfont}
    \begin{lstlisting}
  }{
    \end{lstlisting}
  }
}

\definecolor{jlkeyword}{HTML}{0000FF}%      % Blue for keywords (function, for, if, etc.)
\definecolor{jlfunctions}{HTML}{8B008B}%    % Dark magenta for functions
\definecolor{jlbuiltin}{HTML}{FF4500}%      % Orange red for built-ins
\definecolor{jlmacros}{HTML}{1F7199}%       % Teal for macros
\definecolor{jlliteral}{HTML}{008000}%      % Green for literals (true, false, nothing)
\definecolor{jlstrnum}{HTML}{B22222}%       % Fire brick for strings and numbers
\definecolor{jlcomment}{HTML}{808080}%      % Gray for comments
\definecolor{jlop}{HTML}{2E8B57}%           % Sea green for operators

\usepackage[section]{algorithm}
\usepackage{algorithmicx} 
\usepackage{algpseudocode}

\usepackage{xspace}

\usepackage[figuresright]{rotating}
\usepackage{csvsimple}
\usepackage{makecell}

\usepackage{array}
\newcolumntype{H}{>{\setbox0=\hbox\bgroup}c<{\egroup}@{}}

\usepackage[scaled=0.9]{FiraMono}

\usepackage{booktabs}

\usepackage{amsthm}

\newtheorem{proposition}{Proposition}
\newtheorem{theorem}{Theorem}
\newtheorem{definition}{Definition}
\newtheorem{lemma}{Lemma}
\newtheorem{corollary}{Corollary}

\usepackage{chngcntr}
\usepackage{apptools}
\AtAppendix{

}
 
\newif\ifarxiv
\arxivtrue

\title{Boscia.jl: A review and tutorial}

\author{\name Wenjie Xiao \email \href{mailto:xiao@zib.de}{xiao@zib.de} \\
\addr Technische Universit\"at Berlin, Germany\\
\addr Zuse Institute Berlin, Germany
\AND
\name Deborah Hendrych \email \href{mailto:hendrych@zib.de}{hendrych@zib.de}\\
\addr Technische Universit\"at Berlin, Germany\\
\addr Zuse Institute Berlin, Germany
\AND
\name Mathieu Besançon \email \href{mailto:mathieu.besancon@inria.fr}{mathieu.besancon@inria.fr} \\
\addr Univ.~Grenoble Alpes, Inria, CNRS, LIG, Grenoble, France \\
\addr Zuse Institute Berlin, Germany
\AND
\name Sebastian Pokutta \email \href{mailto:pokutta@zib.de}{pokutta@zib.de} \\
\addr Technische Universit\"at Berlin, Germany\\
\addr Zuse Institute Berlin, Germany
}

\DeclareMathOperator*{\argmin}{argmin}
\DeclareMathOperator{\Tr}{Tr}
\DeclareMathOperator{\rank}{rank}
\DeclareMathOperator{\diag}{diag}

\DeclareMathOperator{\conv}{conv}

\newcommand{\innp}[2]{\left\langle #1, #2 \right\rangle}

\newcommand{\vx}{\mathbf{x}}
\newcommand{\vvv}{\mathbf{v}}

\newcommand{\vy}{\mathbf{y}}

\newcommand{\vc}{\mathbf{c}}
\newcommand{\vl}{\mathbf{l}}
\newcommand{\vu}{\mathbf{u}}
\newcommand{\va}{\mathbf{a}}
\newcommand{\vs}{\mathbf{s}}

\newcommand{\R}{ {\mathbb R} } % defines a short cut for the symbol of real numbers
\newcommand{\Z}{ {\mathbb Z} } % define symbol integers
\usepackage{jmlr2e}
\usepackage[noabbrev,capitalize,nameinlink]{cleveref}
\usepackage{natbib}
\setcitestyle{authoryear,round,citesep={;},aysep={,},yysep={;}}
\usepackage[vvarbb]{newtxmath}
\usepackage{multirow}

\usepackage{todonotes}
\definecolor{brandeisblue}{rgb}{0.0, 0.44, 1.0}
\newif\ifshowrevisions
\showrevisionsfalse

\newcommand{\wenjieplot}[2][]{%
  \begingroup
  \setlength{\fboxsep}{2pt}%
  \setlength{\fboxrule}{0.65pt}%
  \ifshowrevisions
    \fcolorbox{brandeisblue}{white}{\includegraphics[#1]{#2}}%
  \else
    \includegraphics[#1]{#2}%
  \fi
  \endgroup
}

\newcommand{\wenjietable}[1]{%
  \begingroup
  \setlength{\fboxsep}{4pt}%
  \setlength{\fboxrule}{0.65pt}%
  \ifshowrevisions
    \fcolorbox{brandeisblue}{white}{#1}%
  \else
    #1%
  \fi
  \endgroup
}
\newcommand{\wenjiefigure}[1]{%
  \begingroup
  \setlength{\fboxsep}{4pt}%
  \setlength{\fboxrule}{0.65pt}%
  \ifshowrevisions
    \fcolorbox{brandeisblue}{white}{#1}%
  \else
    #1%
  \fi
  \endgroup
}

\newcommand{\revision}[1]{\ifshowrevisions{\color{brandeisblue} #1}\else#1\fi}

\definecolor{revorange}{rgb}{0.0, 0.44, 1.0}
\newif\ifshowrevII
\showrevIIfalse
\newcommand{\revII}[1]{\ifshowrevII{\color{revorange}#1}\else#1\fi}

\begin{document}

\maketitle

\begin{abstract}
    Mixed-integer nonlinear optimization (MINLP) comprises a large class of problems that are challenging to solve
    and exhibit a wide range of structures.
    The Boscia framework \revision{\citep{hendrych2023convex}} focuses on convex MINLP where the nonlinearity appears in the objective only.
    This paper provides an overview of the framework\revision{, showcases extensions post-publication} and practical examples to illustrate its use and customizability.
    One key aspect is the integration and exploitation of Frank-Wolfe methods as continuous solvers within a branch-and-bound framework,
    enabling inexact node processing, warm-starting and explicit use of combinatorial structure among others.
    Three examples illustrate its flexibility, the user control over the optimization process and the benefit of oracle-based access to 
    the objective and its gradient.
    \revision{Additionally, ablation studies are performed on the three examples to investigate the performance impact of the different features and customizations.}
    The aim of this tutorial is to provide readers with an understanding of the main principles of the framework.
\end{abstract}

\section{Introduction}\label{sec:Intro}

Mixed-integer nonlinear programming (MINLP) represents a challenging and wide class of optimization problems, 
combining the combinatorial complexity of integer variables with the nonlinearity in the constraints and the objective function.

The two main approaches to solving (convex) MINLPs are outer approximation (OA) schemes and branch-and-bound (B\&B) 
methods \revision{(see \citet{kronqvist2019review} for a recent review of solution methods)}. 
The former solves a sequence of linear approximations of the original problem, using the gradients of 
the nonlinear constraints and objective function to generate linear cuts to exclude points that are infeasible for the
original problem.
The latter solves a sequence of nonlinear problems (NLP) by dividing the feasible region into smaller subproblems with respect 
to the integer variables. 
Usually, the number of subproblems is exponential in the number of integer variables, and the approach relies
on the subproblems providing good lower bounds.
We focus on convex problems throughout this paper.

When the nonlinearity appears only in the objective and the constraints remain linear, traditional outer-approximation methods may face computational challenges, particularly in higher dimensions. 
\revision{
In particular, this can occur if the constraints encode combinatorial structures that are distorted by the added cuts \citep{hendrych2023convex,2023_HendrychBesanconPokutta_Optimalexperimentdesign}. Standard linearized OA cuts may fail to fully exploit the combinatorial structure of the integer-feasible set and may therefore provide weak lower bounds, as studied in \citet{wei2023outer}\revII{.}
}

The Boscia framework \citep{hendrych2023convex} introduces an approach to solving MINLPs with convex 
objectives and linear constraints by combining \revision{a branch-and-bound algorithm with executions} of Frank-Wolfe (Conditional Gradient) 
\citep{frank1956algorithm,levitin1966constrained} for the continuous subproblems.
It exploits many of the properties of the Frank-Wolfe algorithm to make the B\&B process more efficient,
\revision{which we detail in \cref{sec:Framework}}.
Like Frank-Wolfe, the framework requires only oracle access to the objective function and its gradient, 
as well as to the feasible region in the form of a Linear Minimization Oracle (LMO).
Thus, the framework offers several advantages over traditional MINLP approaches. 
It can exploit combinatorial structure in the linear constraints via the LMO, \revision{leveraging specialized 
algorithms that do not require the algebraic description of the feasible region.
The oracle-based access to the objective means that no closed-form or symbolic representation is necessary.}

\revision{Since its initial release in 2022, the capabilities of the package have been significantly expanded. 
These extensions include support for different Frank-Wolfe variants, new branching strategies, and an extended LMO interface 
supporting MIP solvers and combinatorial algorithms.
Callbacks throughout the stages of the solve allow for more user control and customized behavior, 
and support was added for objectives that are not well-defined over the full feasible region.
}

This paper reviews the Boscia framework—its theoretical foundations, design choices, implementation paradigms, 
and applications—and serves as a practical tutorial to build user intuition and guide implementations.
\revision{In addition, we present new empirical and implementation-focused contributions beyond \citet{hendrych2023convex}. 
\begin{itemize}
    \item We analyze and document the new developments post-publication such as additional branching strategies, broader LMO paradigms, richer callback-based customization, and support for objectives that are not well-defined on the full feasible region.
    \item We perform ablation studies on three different examples to quantify the performance impact of these design choices.
    \item We provide practical guidelines for advanced and customized use cases of Boscia illustrated on three examples, each targeting different aspects of the framework.
\end{itemize}
}

\noindent
\revision{
    Some popular open-source MINLP solvers include SCIP \citep{BolusaniEtal2024OO}, which utilizes polyhedral approximations for 
    nonlinear terms and spatial branching in B\&B for nonconvex problems; BONMIN \citep{bonami2008algorithmic}, which employs 
    both B\&B and outer approximation (and mixtures thereof); and SHOT \citep{lundell2022polyhedral,lundell2022supporting}, 
    which implements outer approximation schemes.
    Popular commercial solvers include Gurobi \citep{gurobi} and Xpress \citep{Xpress}, both also implementing polyhedral approximations.
    }

We begin by examining the framework's architecture and the integration of Frank-Wolfe methods within branch-and-bound 
(\cref{sec:Framework}).
We then present detailed examples demonstrating different implementation approaches, ranging from network 
design problems to graph isomorphism and optimal experiment design problems (\cref{sec:Examples_Results}).
\revision{Each example section includes a computational study to examine the effect of relevant features and customized behavior. 
We conclude the section with practical guidelines for using the framework effectively, derived from the ablation studies.}
Finally, we discuss best practices, current limitations, and future research directions (\cref{sec:Discussion}).

\subsection*{Notation}\label{sec:Notation}

Throughout this work, we assume familiarity with basic optimization concepts (constraints, gradients, convexity) and use the following notation: 
$L$-smooth functions are those whose gradients are Lipschitz continuous with constant $L$ over a given convex compact set, 
and we denote the inner product between vectors $\vx$ and $\vy$ as $\innp{\vx}{\vy}$.
Matrices are denoted by uppercase letters, vectors are in lowercase and bold, scalars in lowercase letters.
Sets are denoted by calligraphic letters.
The convex hull of a set $\mathcal{X}$ is denoted by $\conv(\mathcal{X})$.
\revII{For a vector $\vx$, $\diag(\vx)$ denotes the diagonal matrix with the entries of $\vx$ on its diagonal.}

\section{The Boscia framework}\label{sec:Framework}

The Boscia framework introduced in \citet{hendrych2023convex} solves mixed-integer nonlinear problems (MINLP) with convex objectives 
and linear constraints of the form: 

\begin{align}
    \min_{\vx} &\; f(\vx) \label{eq:minlp} \\
    \text{s.t.} &\; \vx \in \mathcal{X} \nonumber \\
& \; x_i \in \Z \, \forall i \in I \nonumber
\end{align}
\revision{where $\mathcal{X} \subseteq \R^n$ denotes the continuous linear feasibility set. Let $\Z_I = \{\vx \in \R^n \mid x_i \in \mathbb{Z}\;\forall i \in I\}$, $I\subseteq \{1,\ldots,n\}$ denote the set of points satisfying the integrality constraints.}
Further, we will refer to $\vx\in\mathcal{X} \cap \Z_I$ as \emph{integer-feasible}.

The framework employs a branch-and-bound approach with Frank-Wolfe (FW) (alternatively called Conditional Gradient (CG)) \citep{frank1956algorithm,levitin1966constrained} 
methods as the solver for the continuous nonlinear subproblems.
Frank-Wolfe is a first-order algorithm for convex constrained optimization problems that assumes a differentiable and $L$-smooth function $f$.
In the standard version, Frank-Wolfe solves a linear minimization problem in each iteration, using the gradient of the current iterate $\vx_t$ as cost function, 
via a Linear Minimization Oracle (LMO).
\begin{align*}
    \vvv &\;\leftarrow \argmin_{\vy \in \mathcal{X}} \innp{\nabla f(\vx_t)}{\vy}
\end{align*}
The returned extreme point of the (continuous) feasible region $\mathcal{X}$ is then used to update the current iterate $\vx_t$ via a convex combination:
\begin{align*}
    \vx_{t+1} &=\vx_t + \gamma_t (\vvv - \vx_t)
\end{align*}
where $\gamma_t$ is the step size, either fixed depending on the iteration or computed via a line search.
The algorithm continues until the \emph{FW gap} defined below is sufficiently small:
\begin{align*}
    g(\vx) = \max_{\vvv \in \mathcal{X}} \innp{\nabla f(\vx)}{\vx - \vvv}.
\end{align*}
Two things are noteworthy. First, the FW gap upper-bounds the primal gap of the continuous problem,
\begin{align}
    \label{eq:lower_bound}
    f(\vx) -f^* \leq g(\vx) \quad \forall \vx \in \mathcal{X}
\end{align}
where $f^*$ is the continuous optimal value.
This follows from convexity of $f$ and optimality of the computed vertex $\vvv$ for the linear subproblem. 
Consequently, $f(\vx) - g(\vx)$ is a valid lower bound on the optimal objective value.
The FW gap is also driven to zero at optimality and converges along with the primal gap. 
A comprehensive proof can be found in \citet{jaggi2013revisiting}.
The point $\vvv$ is precisely the one computed by the LMO and hence, the FW gap is practically computed for free as a by-product of all iterations.

There are several variants of the Frank-Wolfe algorithm, for an overview see \citet{braun2022conditional}.
Most of the variants currently supported by our framework are \emph{active-set-based}. 
The \emph{active set} encodes the convex combination of extreme points (with their corresponding weights) forming the current iterate.
Actively storing this active set enables the variants to take steps within the active set, thereby keeping the support
of the solution small.
% Algorithm 1: Corrective Frank-Wolfe
\begin{algorithm}[H]
    \caption{Corrective Frank-Wolfe (CFW) \citep{halbey2025efficient}}
    \label{alg:cfw}
    \begin{algorithmic}[1]
    \Require convex, smooth function $f$, start point $\vx^0 \in \mathcal{V}(\mathcal{X})$ (vertex of $\mathcal{X}$).
    \State $\mathcal{S}_0\gets \{\vx_0\}$ \Comment{active set}
    \For{$t = 0$ to $T-1$}
      \State $\va_t \gets \arg\max_{\vvv \in \mathcal{S}_t} \langle \nabla f(\vx_t), \vvv \rangle$ \Comment{away vertex} \label{alg:cfw:a_t}
      \State $\vs_t \gets \arg\min_{\vvv \in \mathcal{S}_t} \langle \nabla f(\vx_t), \vvv \rangle$ \Comment{local FW} \label{alg:cfw:s_t}
      \State $\vvv_t \gets \arg\min_{\vvv \in \mathcal{V}(\mathcal{X})} \langle \nabla f(\vx_t), \vvv \rangle$ \Comment{global FW} \label{alg:cfw:v_t}
      \If{$\langle \nabla f(\vx_t), \va_t - \vs_t \rangle \geq \langle \nabla f(\vx_t), \vx_t - \vvv_t \rangle$} \label{alg:cfw:local_gap}
        \State $\vx_{t+1}, \mathcal{S}_{t+1} \gets \textsc{CorrectiveStep}(\mathcal{S}_t, \vx_t, \va_t, \vs_t)$
      \Else
        \State $\gamma_t \gets \arg\min_{\gamma \in [0,1]} f(\vx_t - \gamma(\vx_t - \vvv_t))$
        \State $\vx_{t+1} \gets \vx_t - \gamma_t (\vx_t - \vvv_t)$
        \State $\mathcal{S}_{t+1} \gets \mathcal{S}_t \cup \{\vvv_t\}$
      \EndIf
    \EndFor
    \end{algorithmic}
    \end{algorithm}
  \vspace{-0.7em}
  % Algorithm 2: CorrectiveStep
  \begin{algorithm}[H]
    \caption{Corrective Step($\mathcal{S}, \vx, \va, \vs$) \citep{halbey2025efficient}}
    \label{alg:corrective-step}
    \begin{algorithmic}[1]
    \Require $\mathcal{S} \subset \mathcal{X}, \; \vx, \va, \vs \in \mathcal{X}$
    \Ensure $\mathcal{S}' \subseteq \mathcal{S}, \; \vx' \in \conv(\mathcal{S}')$ satisfying
      % first alternative (drop step)
      \Statex \(\displaystyle f(\vx') \le f(\vx) \quad\text{and}\quad \mathcal{S}' \subsetneq \mathcal{S} \quad \textbf{{or}}  \)
      \hfill \(\triangleright\) \text{drop step}
      % second alternative (descent step)
      \Statex \(\displaystyle
      f(\vx) - f(\vx') \;\ge\;
      \frac{\langle \nabla f(\vx),\, \va-\vs\rangle^2}{2 L D^2} \)
      \hfill \(\triangleright\) \text{descent step}
    \end{algorithmic}
    \end{algorithm}
All these variants can be interpreted as different cases of the Corrective Frank-Wolfe algorithm (CFW) from \citet{halbey2025efficient}, 
shown in Algorithm \ref{alg:cfw} and Algorithm \ref{alg:corrective-step}.
The active-set-based FW variants differ primarily in the operations performed in the corrective step.
Note that $D$ in Algorithm \ref{alg:corrective-step} denotes the diameter of the set $\mathcal{X}$; 
$L$ is the Lipschitz constant.
In the \emph{Blended Pairwise Conditional Gradient (BPCG)}, for example, the corrective step consists of shifting weight from the away vertex
$\va_t$ to the local FW vertex $\vs_t$.
All of the active-set-based variants can be lazified, meaning the LMO is called only when no further local progress can be made, 
see \cite{braun2017lazifying} for details.

In contrast to classic B\&B approaches, the subproblems solved at each node are not the continuous relaxations of \revision{original problem, that is, }just dropping the integrality constraints. 
Instead, the \revision{objective} function $f$ is optimized over the convex hull of integer-feasible points\revision{, denoted by} $\conv(\mathcal{X})$.
That is, the extreme points returned by the LMO are integer-feasible points, which is obtained by propagating the integrality constraints to the LMO.
So, the LMO encodes \revision{the following} mixed-integer linear problem (MILP) \revision{at a  given node with bounds $[\vl, \vu]$:
\begin{align*}
    \min_{\vy} &\; \innp{\vc}{\vy} \\
    \text{s.t. } &\; \vy \in \mathcal{X} \cap \Z_I \cap [\vl,\vu]
\end{align*}}
which has important computational implications.
\revision{Note that if the feasible region is uni-modular, then $\conv(\mathcal{X}\cap \Z_I)$ and the continuous relaxation of $\mathcal{X}$ are identical.
While the additional integrality constraints make the LMO, and by extension, the Frank-Wolfe algorithm more expensive, it leads to an inbuilt heuristic as any extreme point is a potential solution. Moreover, many combinatorial structures admit to efficient algorithms, for example the Birkhoff polytope and the Matching polytope.}
A schematic of the framework is shown in \cref{fig:schematic}.

Following the Frank-Wolfe model, we assume an oracle for the objective function $f$ and its gradient $\nabla f$.
Our framework requires additional structure for the LMO since it needs to handle the integrality constraints, as detailed in \cref{sec:Framework:BLMO}.

A key consequence is that the LMO call can be relatively expensive, particularly for generic polytopes.
We mitigate this computational burden through several mechanisms. First, Frank-Wolfe is error-adaptive, meaning the computational 
load increases with precision requirements.
We exploit this by using looser tolerances near the root and tightening them deeper in the tree.
This error-adaptivity is quite rare in MINLP solvers which typically require solving subproblems almost to optimality at every node.
By \eqref{eq:lower_bound}, Frank-Wolfe always provides a valid lower bound on the continuous problem, and 
thus, a valid lower bound for the tree.

\begin{wrapfigure}{r}{0.5\textwidth}
    \centering
    \begin{tikzpicture}[
        node distance=1cm and 1cm,
        every node/.style={font=\small},
        >=Stealth,
        box/.style={rectangle, draw, minimum width=3cm, minimum height=1cm, align=center},
        rhombus/.style={shape=diamond, draw, aspect=1.3, minimum width=2.0cm, minimum height=0.8cm, align=center, inner sep=2pt},
        ellipsoid/.style={shape=ellipse, draw, minimum width=2cm, minimum height=1.5cm, align=center},
        arrow/.style={->, thick, bend angle=25}
    ]
    
    % Nodes
    \node[box] (top) {\textbf{B\&B}\\ \\$\begin{array}{rl}
        \min_\vx & f(\vx) \\
        \text{s.t.} & \vx \in\mathcal{X} \revision{ \cap \Z_I}
    \end{array}$};
    \node[ellipsoid, below=of top] (middle) {\textbf{FW}\\ \\$\begin{array}{rl}
        \min_\vx & f(\vx) \\
        \text{s.t.} & \vx \in\conv(\mathcal{X} \revision{\cap \Z_I})
    \end{array}$};
    \node[rhombus, below=of middle] (bottom) {\textbf{LMO}\\ \\$\begin{array}{rl}
        \min_\vy & \innp{\nabla f(\vx)}{\vy} \\
        \text{s.t.} & \vy \in \mathcal{X} \revision{\cap \Z_I}
    \end{array}$};
    
    % Arrows between top and middle (rhombus)
    \draw[arrow, bend right] (top.south west) to node[left]{node} (middle.north west);
    \draw[arrow, bend right] (middle.north east) to node[right]{bound} (top.south east);
    
    % Arrows between middle and bottom (ellipsoid)
    \draw[arrow, bend right] (middle.south west) to node[left]{gradient} (bottom.north west);
    \draw[arrow, bend right] (bottom.north east) to node[right]{direction} (middle.south east);
    
    \end{tikzpicture}
    \caption{The schematic of the algorithm including the optimization problems solved at the different layers.}
    \label{fig:schematic}
\end{wrapfigure}

Second, as stated previously, the Frank-Wolfe variants employed are active-set-based.
This information can be utilized during branching by splitting the active set into two parts, one for the left and one for the right child.
Note that since we assume that our vertices are always integer-feasible, none of the new active sets can be empty.
This enables warm-starting the child nodes, reducing the number of iterations needed to convergence to the desired tolerance.
Third, we utilize the lazification mechanism initially proposed in \citet{braun2017lazifying}, meaning we can use vertices that are not minimizers of the linear subproblems as long as they provide progress.
We also add tree-level lazification by maintaining a shadow set that keeps track of discarded vertices and stores them in a pool passed to child nodes.
While these vertices are not useful for the current problem since they were discarded, they may be useful for the child nodes and can thus avoid being recomputed through the LMO.

Adding the integrality constraints to the LMO subproblems might at first glance seem expensive, since it will often mean solving NP-hard problems. 
However, by doing so, we obtain integer-feasible points from the root node, and thus obtain an upper bound for the tree. 
This leads to a significantly smaller search tree as many nodes can be pruned early.

The full algorithm can be found in \citet[Algorithm~2.1]{hendrych2023convex} as well as performance results.
The framework is implemented in \texttt{Julia} and is available as \package{} \citep{BosciaGitHub}.

\subsection{Frank-Wolfe variants}\label{sec:Framework:FW}

In the following, we give a brief overview of the different Frank-Wolfe variants supported by \package{}
which are implemented in the \texttt{FrankWolfe.jl} package \citep{besanccon2025improved,besanccon2022frankwolfe,FrankWolfeGitHub}.
First, we have the \emph{Standard Frank-Wolfe}, albeit with an active set which can be used for lazification and 
also for warm-starting the child nodes.

The \emph{Away-Frank-Wolfe (AFW)} variant \citep{wolfe1970convergence} was the first variant to be proposed.
Its main idea is to mitigate the zig-zagging behavior displayed by standard Frank-Wolfe, if the optimal solution is on a face of the feasible region.
By introducing away steps, the algorithm is able to move away from a suboptimal vertex or even to drop it entirely from the convex combination of the iterate.

The \emph{Blended Conditional Gradient (BCG)} \citep{braun2019blended} variant's corrective step is a descent step over the 
convex hull of the active set.
If the descent step is computed to optimality, we have the Fully-Corrective-Frank-Wolfe (FCFW), first introduced in \citet{holloway1974extension}.
Note that while FCFW makes a lot of progress per iteration, the optimization over the convex hull of the active set is expensive in wall clock time.

Extending the idea of away steps, the \emph{Pairwise Frank-Wolfe (PFW)} variant’s corrective step moves between the away 
vertex and the global FW vertex, i.e., the LMO solution.
The theoretical convergence rate is not very good due to swap steps—steps during which the away vertex is discarded. 
In practice, however, this variant has shown to be quite efficient.

The \emph{Blended Pairwise Conditional Gradient (BPCG)} variant \citep{braun2019blended}, already discussed above, performs pairwise steps like PFW but only between vertices of the active set.
If this pairwise step does not provide enough progress, a standard FW step is performed. 
This variant has good theoretical convergence rate and is efficient in practice, hence it is the default variant in \package{}.

\revision{Lastly, we consider the active-set-free \emph{Decomposition-Invariant Conditional Gradient (DICG)} method \citep{garber2016linear}.}
The motivation for this variant is two drawbacks of the active-set-based variants: 
First, the convex combination is, in general, not unique and the performance of the methods can vary greatly depending on the quality of the combination.
Second, for large dimensional problems, we might run into the storage issues of the active-set-based methods.
On the one hand, the vertices to be stored are large and on the other hand, the active set tends to grow with the dimension of the problem.
The idea of DICG is to compute the away vertex by restricting the LMO to the minimal face containing the current iterate. 
The rationale is that this yields the worst possible vertex over all convex combinations.
This is referred to as the \emph{in-face LMO} call.
The minimal face is defined as the smallest face of $\mathcal{X}$ containing $\vx$.
The weight is then shifted from this away vertex to the global FW vertex. 
Note that this variant still implicitly uses a convex combination.
DICG is effective on structured polytopes, like simplices, the Birkhoff polytope and hypercubes. 
Its main drawbacks are threefold: (1) identifying the minimal face may be expensive; (2) computing the in-face away vertex 
requires another LMO call (usually in smaller dimension); and (3) the maximum step size cannot be read off the active-set 
weights and must be computed explicitly.
For an arbitrary polytope, active-set-based methods are therefore preferable.
\revision{Note that Boscia also supports the \emph{Blended Decomposition-Invariant Conditional Gradient (BDICG)} variant \citep{besanccon2025improved}, which includes the logic of the blended step from BPCG in DICG.}

The user can also easily add support for their own FW variant by following the same interface used to integrate variants from the \texttt{FrankWolfe.jl} package.

\subsection{The LMO with integer constraints}\label{sec:Framework:BLMO}

In Boscia, the LMO interface from \texttt{FrankWolfe.jl} has to be extended to support the propagation of the integrality constraints to the LMO
\begin{align}
    \vvv &\;\leftarrow \argmin_{\vy \in \mathcal{X} \cap \Z_I \cap [\vl,\vu]} \innp{\vc}{\vy} \label{eq:blmo}
\end{align}
where $\mathcal{X}$ is the \revision{continuous} feasible region of the problem, $\Z_I$ the set of integer points, and $[\vl,\vu]$ are 
the local node bounds on the integer variables.
\revision{Hence, the returned point $\vvv$ is always integer-feasible and an extreme point of $\conv(\mathcal{X}\cap\Z_I)$. 
This provides the B\&B with upper bounds and possible solutions from the first Frank-Wolfe iteration at the root node. 
We have observed that often the solution is found early as a vertex of the low-depth subproblem and the B\&B spends most of its solving time proving optimality.}

The key \revision{implementation} innovation lies in the dynamic construction of LMOs for each branch-and-bound node, where integrality constraints are propagated through 
bound management rather than maintaining separate copies of the LMO for each node. 
This design choice significantly reduces memory overhead while maintaining computational efficiency.
On the other hand, it requires more functionality from the LMO, namely the ability to read, set, add and delete bounds.
Additionally, there are feasibility checks and performance logging functions that while optional can be useful to implement.

The framework offers three distinct implementation\revision{/modelling} pathways.
\revision{The most straightforward way are generic MIP solvers, like \texttt{SCIP} \citep{BolusaniEtal2024OO} or \texttt{HiGHS} \citep{huangfu2018parallelizing}, via Julia's
modelling tool \texttt{JuMP.jl} \citep{Lubin2023} or its backend \texttt{MathOptInterface.jl} (MOI, \citet{mathoptinterface}). 
Any solver supporting MOI can be used; provided it can solve Problem \eqref{eq:blmo}.}

If the problem in \eqref{eq:blmo} \revision{admits to an efficient combinatorial algorithm}, the bound 
management \revision{can be} handled by the framework \revision{itself} and only the computation of the bounded extreme point and 
a feasibility check for \revision{linear} constraints have to be implemented.
\revision{The package already contains several such simple LMOs, for example for the cube, the box and the simplex.
The solver then only requires the set of integer variables and, if applicable, their global bounds.}
For an example, we refer the reader to \cref{sec:Examples:NetworkDesign}.

\revision{If the bound management is performance-critical, because, for example, bound changes have direct effects on internal structures, it is advisable 
to implement the full LMO interface. 
One such example is the Birkhoff polytope as shown in the example in \cref{sec:Examples:GraphIsomorphism}. 
The corresponding LMO can be found in the \texttt{CombinatorialLinearOracles.jl} package \citep{CLOraclesGitHub} which also includes LMOs for other combinatorial
polytopes such as the matching polytope.}

Each approach balances ease of implementation with computational performance, allowing users to choose the most appropriate method based 
on their specific problem structure and computational requirements.

\subsection{Optional settings}\label{sec:Framework:Settings}

The framework is designed with flexibility and user control as primary objectives, providing extensive configurability throughout the solving process.
Node and time limits are standard branch-and-bound parameters. 
As a traversal strategy, the node with the smallest lower bound is selected.
The framework includes a number of branching strategies, varying in complexity and performance:
\begin{itemize}
    \item most-infeasible branching\revision{, also known as most-fractional branching,} which is simple and relatively effective for many problems \revision{\citep{bonami2013branching}},
    \item strong branching is costly since a few FW iterations are performed for each potential branching candidate \revision{\citep{applegate1998solution}},
    \item pseudo-cost branching strategies adapted to the nonlinear case \revision{\citep{benichou1971experiments,bonami2013branching,belotti2009branching}},
    \item gradient-based branching \revision{(\citet{belotti2009branching} use gradient information for scoring)},
    \item and hierarchy branching where the branching strategies are applied successively in a user-specified order \revision{\citep{achterberg2009hybrid}}.
\end{itemize}

The framework also incorporates a number of callback mechanisms.
The B\&B callback is called right before the next node is evaluated and is utilized, for example, for logging.
This enables users to monitor progress, extract the incumbent solution, and stop the algorithm early.
Note that any user provided callback will be called before the internal callback mechanism.
Additionally, a specialized branch callback, called during branching, allows users to selectively prevent the creation of child nodes, 
providing fine-grained control over the search tree exploration, for an example 
see \cref{sec:Examples:GraphIsomorphism}.
In \citet{2025_MexiEtAl_Frankwolfeheuristic_2508-01299}, the various callbacks and node and time limits are used to implement a heuristic framework for quadratic problems based on Boscia.

The framework includes many Frank-Wolfe specific settings that allow users to customize the optimization algorithm to their needs:
\begin{itemize}
    \item Frank-Wolfe variants,
    \item line search (all of line search methods supported by \texttt{FrankWolfe.jl}),
    \item maximum number of iterations and time limit ,
    \item use of lazification and sparsity control.
\end{itemize}
The standard line search is the Secant line search \citep{hendrych2025secant}, which was shown to be particularly efficient on 
quadratic and self-concordant functions \citep{sun2019generalized}. 
\revII{
It seeks a step size $\gamma$ satisfying
$\langle \nabla f(x_t+\gamma d_t), d_t\rangle = 0$, i.e., a stationary point of the objective along the search direction $d_t$.
For quadratic and self-concordant objectives, this allows the optimal step size to be identified efficiently.
The framework further supports the Agnostic step-size rule \citep{jaggi2013revisiting} and the Golden Ratio line search \citep{kiefer1953sequential,Avriel01101966}.
The former uses $\gamma_t=\ell/(\ell+t)$ for a prescribed parameter $\ell$.
The latter applies golden-section search to the function $\phi(\gamma)=f(x_t+\gamma d_t)$ over the admissible step-size interval, progressively contracting an interval containing a minimizer without requiring derivative or smoothness information.

A second class of strategies relies on the smoothness constant $L$ or a local estimate $M$.
The Short Step line search chooses
$\gamma_t=\min\{g_t/(L\|d_t\|^2),\gamma_t^{\max}\}$,
where $g_t=\langle-\nabla f(x_t),d_t\rangle$.
When $L$ is unknown, both the Backtracking and Adaptive line searches replace it by a local estimate $M$ and use the corresponding short-step formula.
The Backtracking line search \citep{pedregosa2020linearly} increases $M$ until the sufficient-decrease condition
$f(x_t+\gamma_t d_t)\leq Q_t(\gamma_t,M)$
is satisfied, where
$Q_t(\gamma,M)=f(x_t)+\gamma\langle\nabla f(x_t),d_t\rangle+\frac{M}{2}\gamma^2\|d_t\|^2$
is the quadratic smoothness model along $d_t$.
In contrast, the Adaptive line search \citep{pokutta2024frank} uses the same candidate step size, but accepts the estimate $M$ whenever
$\langle\nabla f(x_{t+1}),x_t-v_t\rangle\geq 0$.
}

Tolerance settings provide another layer of control, including relative and absolute tolerances for the branch-and-bound process which 
control termination of the algorithm.
Available tolerances include:
\begin{itemize}
    \item absolute tolerance (default $10^{-6}$) and relative tolerance (default 1\%),
    \item FW gap decay factor (for adaptive tolerance tightening over the tree, default 0.8),
    \item start FW epsilon and minimum FW epsilon (defaults $10^{-2}$ and $10^{-6}$, respectively),
    \item min lower bound (unset by default).
\end{itemize}
Since FW provides a valid lower bound in each iteration, there is the option to prematurely stop the evaluation of a
node if enough other open nodes have a better initial lower bound. 

A post-processing procedure is in place for mixed-integer problems.
It fixes the integer variables to the values of the best solution and runs the chosen FW variant only for the continuous variables.

Heuristic strategies are implemented with probabilistic activation, where each heuristic has an associated probability that determines whether it should be applied at each node.
Due to the probabilistic activation, computationally expensive heuristics can be applied only occasionally.
The framework includes several built-in heuristics:
\begin{itemize}
    \item simple rounding heuristic
    \item probability rounding heuristic
    \item follow-the-gradient heuristic 
    \item hyperplane-aware heuristics for simplex-like feasible regions
\end{itemize} with simple rounding activated by default while others remain optional.
The heuristic interface enables easy integration of custom heuristics.
An optional solution callback is triggered whenever a new solution is added to the search tree.

There are gradient-based tightening strategies applicable either globally or locally.
Additionally, strong convexity and sharpness can also be exploited to tighten bounds.

Finally, domain settings address scenarios where the function is not well-defined over the whole feasible region, 
providing mechanisms to handle issues arising from such constraints, as demonstrated in the optimal experimental 
design problem example, see \cref{sec:Examples:OEDP}.

\section{Examples and computational results}\label{sec:Examples_Results}

In the following, we present three examples that highlight different customization and workflows with the 
Boscia framework. 
\revision{Furthermore, we perform ablation studies to investigate the impact of different Boscia features and customizations on the performance of the framework.}
The first is a network design problem which showcases two different ways of modeling the feasible region and is a 
good baseline example.
\revision{In our experiments, we compare two different formulations of the problem with different LMO models. 
Additionally, we compare the performance of different Frank-Wolfe variants, line-search 
methods, and branching strategies.}
The second example is the graph isomorphism problem where we demonstrate how to customize the solving process 
via the available callback mechanism and how to implement a fully user-managed LMO.
\revision{In the experiment setup, we focus on the performance impact of the different LMO implementations, comparing the MIP-solver-based LMO, 
a solver-managed and self-managed LMO.}
The last example deals with the case of the objective function not being well-defined over the whole feasible
region based on the optimal experiment design problem. 
\revision{We then report the impact of different heuristic strategies, Frank-Wolfe tolerance settings, and lazification with 
and without the shadow set.}
The source code for all examples can be found on the GitHub repository \citep{BosciaGitHub} and in the
\href{https://zib-iol.github.io/Boscia.jl/stable/}{documentation}~\citep{BosciaDocumentation}.

\paragraph{\revision{Environment}}

\revision{
All experiments are implemented in Julia. We use \texttt{Boscia.jl} v0.2.11 and 
\texttt{FrankWolfe.jl} v0.6.4 for the Frank-Wolfe node relaxations. Solver interfaces are provided 
through \texttt{MathOptInterface.jl} v1.11.0 and \texttt{MathProgBase.jl} v1.47.0. Mixed-integer 
linear subproblems are solved with \texttt{SCIP.jl} v0.12.7 or \texttt{HiGHS.jl} v1.20.1.
}

\revision{
The experiments are run on an OpenStack virtual compute node with 128 Intel Xeon Cascadelake 
CPU cores and 32 GB of RAM. Unless stated otherwise, all experiments use this hardware 
configuration and the package versions listed above. Time limits and termination criteria are 
specified for each problem class separately.
}

\subsection{Network design problem - simple LMO and MOI LMO}\label{sec:Examples:NetworkDesign}

Suppose we have given a network $G= (\mathcal{V}, \mathcal{E})$, a list of potential edges $\mathcal{R}$ not in $\mathcal{E}$, and a list of 
flow demands between sources $\mathcal{O}\subset \mathcal{V}$ and destinations $\mathcal{Z}\subset \mathcal{V}$. 
The goal is to minimize the design cost of adding new edges from the potential edge set $\mathcal{R}$ and 
the operating cost of the network.
This is a good baseline example because it requires few customizations.
It can be modeled using two different LMOs and thus, we can showcase the two most common approaches to LMOs in \package{}.
The example is based on \citet{2024_SharmaHendrychBesanconPokutta_NetworkdesignMicoFrankwolfe} to which we refer the reader 
for comprehensive computational results. 
Here, the operating cost of the network is modeled as a \emph{traffic assignment problem} (TA) with a congestion effect.
\begin{align}
    \min_{\vx} &\; c(\vx) := \sum_{e \in \mathcal{E}} c_e(x_e) \tag{TA}\label{eq:ta}\\
    \text{s.t.} &\; x_e = \sum_{z \in \mathcal{Z}}x_e^z \;&&\forall e \in \mathcal{E}\nonumber \\
    &\; \vx^z \in \mathcal{X}^z = \begin{cases}
        \sum_{e \in \delta^{+}(i)}x_e^z - \sum_{e \in \delta^{-}(i)}x_e^z = 0, \;\forall i \in \mathcal{V} \ \backslash \ (\mathcal{O} \cup \mathcal{Z})\\
        \sum_{e \in \delta^{+}(i)}x_e^z = d_i^z \;\forall i \in \mathcal{O}\\
        \sum_{e \in \delta^{-}(z)}x_e^z = \sum_{i \in \mathcal{O}} d_i^z\\
        \end{cases}\;&&\hspace{-1mm}\forall z \in \mathcal{Z}.\nonumber
\end{align}
The first constraint simply ensures that the total flow on an edge is the sum of the flows to all destinations.
Observe that we do not assume any capacity constraints on the edges.
At all nodes that are neither sources nor destinations, the flow has to be balanced.
The sum of outgoing flows from a source $i$ to a destination $z$ has to be equal to the demand between $i$ and $z$.
Likewise, the sum of all incoming flows to a destination $z$ has to be equal to the sum of all demands to $z$.
The edge cost functions $c_e$ estimate the travel time and are modeled as:
\[ c_e(x_e)= \alpha_e + \beta_ex_e + \gamma_e x_e^{\rho_e}\]
where $\alpha_e, \beta_e$ and $\gamma_e$ are constants and the exponents $\rho_e > 1$ model the congestion
effect of a network.

For the network design problem, we add binary variables $\vy \in\mathcal{Y}\subseteq \{0,1\}^{|\mathcal{R}|}$ 
that model which edges from $\mathcal{R}$ should be added to the network and a linking constraint per edge $e$ 
in $\mathcal{R}$ forcing the corresponding flow on the edge to zero if the associated binary is zero:
\begin{align}
    \min_{\vy, \vx} \;\; & \mathbf{r}^\intercal \vy + c(\vx) \tag{ND}\label{eq:ndta}\\
    \text{s.t.} \;\; & y_e = 0 \Rightarrow x_e \leq 0 && \forall e \in \mathcal{R} \nonumber\\
    & \vx \in \mathcal{F}\nonumber \\
    & \vy \in \mathcal{Y} \nonumber
\end{align}
where $\mathcal{F}$ encodes the flow constraints from traffic assignment problem \eqref{eq:ta}.
A small example is shown in \cref{fig:traffic_flow}.

\begin{figure}[htbp]
\centering
\subfloat[Initial Traffic Flow]{
\begin{tikzpicture}[>=Stealth, node distance=1cm, font=\tiny, transform shape, scale=0.9]
    % Nodes
    \node[rectangle, draw, fill=white, text=red] (s1) at (0,2) {$S_1$};
    \node[circle, draw, fill=white, right=of s1] (n1) {1};
    \node[circle, draw, fill=white, below=of n1] (n2) {2};
    \node[circle, draw, fill=white, right=of n1] (n3) {3};
    \node[circle, draw, fill=white, right=of n3] (n4) {4};
    \node[circle, draw, fill=white, below=of n4] (n5) {5};
    \node[circle, draw, fill=white, below=of n2, text=blue] (d) {$D$};
    \node[rectangle, draw, fill=white, below=of n3, text=red] (s2) {$S_2$};

    % Flows
    \draw[->, line width=1pt] (s1) -- (n1) node[midway, above] {$1$};
    \draw[->, line width=0.5pt, dashed] (n1) -- (n2) node[midway, left] {$0$};
    \draw[->, line width=1pt] (n1) -- (n3) node[midway, above] {$1$};
    \draw[->, line width=2.0pt] (n3) -- (n4) node[midway, above] {$2$};
    \draw[->, line width=2.0pt] (n4) -- (n5) node[midway, left] {$2$};
    \draw[->, line width=2.0pt] (n5) -- (d) node[pos=0.4, below=0.15cm] {$2$};
    \draw[->, line width=0.5pt] (n2) -- (d) node[midway, left] {$0$};
    \draw[->, line width=1pt] (s2) -- (n3) node[midway, left] {$1$};
\end{tikzpicture}
}
\hspace{2cm}
\subfloat[Optimized Traffic Flow]{
\begin{tikzpicture}[>=Stealth, node distance=1cm, font=\tiny, transform shape, scale=0.9]
    % Nodes
    \node[rectangle, draw, fill=white, text=red] (s1) at (0,2) {$S_1$};
    \node[circle, draw, fill=white, right=of s1] (n1) {1};
    \node[circle, draw, fill=white, below=of n1] (n2) {2};
    \node[circle, draw, fill=white, right=of n1] (n3) {3};
    \node[circle, draw, fill=white, right=of n3] (n4) {4};
    \node[circle, draw, fill=white, below=of n4] (n5) {5};
    \node[circle, draw, fill=white, below=of n2, text=blue] (d) {$D$};
    \node[rectangle, draw, fill=white, below=of n3, text=red] (s2) {$S_2$};

    % Flows
    \draw[->, line width=1pt] (s1) -- (n1) node[midway, above] {$1$};
    \draw[->, line width=1pt] (n1) -- (n2) node[midway, left] {$1$};
    \draw[->, line width=0.5pt] (n1) -- (n3) node[midway, above] {$0$};
    \draw[->, line width=1pt] (n3) -- (n4) node[midway, above] {$1$};
    \draw[->, line width=1pt] (n4) -- (n5) node[midway, left] {$1$};
    \draw[->, line width=1pt] (n5) -- (d) node[pos=0.4, below=0.15cm] {$1$};
    \draw[->, line width=1pt] (n2) -- (d) node[midway, left] {$1$};
    \draw[->, line width=1pt] (s2) -- (n3) node[midway, left] {$1$};
\end{tikzpicture}
}
\caption{A traffic flow optimization example reproduced from \citet{2024_SharmaHendrychBesanconPokutta_NetworkdesignMicoFrankwolfe} showing initial and 
optimized flow distributions. The line thickness represents flow magnitude, with dashed lines indicating optional edges that can be added.}
\label{fig:traffic_flow}
\end{figure}

Problem \eqref{eq:ta} can be solved in \package{} using a MIP solver, e.g.~\texttt{SCIP} 
\citep{BolusaniEtal2024OO}, supporting the \texttt{JuMP} modeling language in Julia. 

Note that an objective does not need to be set; it is supplied during the solving process by the framework.
Once the optimizer is set with \texttt{JuMP} or \texttt{MathOptInterface}, it has to be wrapped in a dedicated 
LMO for \texttt{JuMP} models.
\begin{jllisting}[language=julia, style=jlcodestyle]{}
lmo = FrankWolfe.MathOptLMO(optimizer)
\end{jllisting}
Next is the definition of the function and its gradient. Neither Boscia nor the underlying Frank-Wolfe
methods require a symbolic expression of the cost function or the gradient. 
They can be given as black-box functions.
For the function, it receives the point to evaluate the function on and we expect the objective value as a return value.
For the gradient, it receives the storage array to write the gradient into for effective memory usage 
and the point at which to evaluate the gradient.
Finally, we can call the solver itself:

\noindent
\begin{minipage}{\linewidth}
\begin{jllisting}[language=julia, style=jlcodestyle]{}
settings = Boscia.create_default_settings()
settings.branch_and_bound[:verbose] = true
x, tlmo, result = Boscia.solve(f, grad!, lmo, settings=settings)
\end{jllisting}
\end{minipage}

The \lstinline{solve} call only requires the function, its gradient\footnote{Julia convention: If functions change their input parameters,
an exclamation mark is appended to the function name.} and the LMO.
The settings will be created by default but can be customized as shown above.
The algorithm returns the solution, the LMO wrapped in a \lstinline{TimeTrackingLMO}\footnote{\revision{A wrapper applied by the framework 
which tracks statistics like number of LMO calls and running time of individual LMO calls.}} object and a result dictionary containing logging and trajectory data.

A downside to the model \eqref{eq:ndta} is the presence of the linking constraints $y_e = 0 \Rightarrow x_e \leq 0$.
The larger the network, the more expensive the call to a MIP solver becomes. 
The performance of the framework depends on the computational cost of the LMO, hence relatively cheap LMOs are preferred.
Thus, \citet{2024_SharmaHendrychBesanconPokutta_NetworkdesignMicoFrankwolfe} suggested an alternative
formulation based on a penalty approach.
\begin{align}
    \min_{\vy, \vx} &\; \mathbf{r}^\intercal \vy + c(\vx) + \mu \sum_{z \in \mathcal{Z}} \sum_{e \in \mathcal{R}} \max(x_e^z - M^z y_e, 0)^p  \tag{PB-ND}\label{eq:rndta}\\
    \text{s.t.} 
    &\; \vx \in \mathcal{F}\nonumber\\
    &\; \vy \in \mathcal{Y}\nonumber
\end{align}
The linking constraints are moved to the objective via a penalty term.
As suggested in \citet{2024_SharmaHendrychBesanconPokutta_NetworkdesignMicoFrankwolfe}, we set $\mu = 10^3$ and $p = 1.5$.
Note that the linking constraints are modeled via a big-M formulation for this approach where $M$ is an upper bound
on the possible flow.
Now, the LMOs of the design variables $\vy$ and the flow variables $\vx$ are computed independently. 
To compute the LMO for the flow, we solve the shortest path problem for each pair of source-destination
and add the demand as flow to all edges in that path.
Note that this is valid since we assume no capacity constraints on the edges.
The LMO for the design variables depends on $\mathcal{Y}$, assumed to be the unit cube in this example.

Since problem \eqref{eq:blmo} can be computed efficiently, the bound management can be left to the framework.
\revision{The LMO only has to implement two functions: one that computes an extreme point under the current node bounds, and one 
that checks feasibility with respect to the simple linear constraints. Their interfaces are shown in Appendix~\ref{app:network_design}.}

The LMO can now be wrapped in a \lstinline{ManagedLMO} which will handle the bound management.
\begin{jllisting}[language=julia, style=jlcodestyle]{}
managed_lmo = Boscia.ManagedLMO(lmo, lower_bounds, upper_bounds, 
    int_vars, total_vars)
x, tlmo, result = Boscia.solve(f, grad!, managed_lmo)
\end{jllisting}

\subsubsection*{\revision{Computational results}}

\revision{
To better study how problem size affects the performance of the LMOs, we consider the stochastic variant of the problem from \citet{2024_SharmaHendrychBesanconPokutta_NetworkdesignMicoFrankwolfe}. 
In contrast to the deterministic network design problem in \cref{eq:ndta}, the stochastic variant accounts for uncertain demand and therefore provides a more realistic model when future traffic levels are not known in advance.
Let $\mathcal{S}$ be a finite set of scenarios with probabilities $p_s$. For each scenario 
$s \in \mathcal{S}$, we denote the flow variables by $\vx_s$ and the corresponding traffic-assignment 
feasible region by $\mathcal{F}_s$. The stochastic network design problem is
\begin{align}
    \min_{\vy, \vx} \;\; 
    & \mathbf{r}^\intercal \vy + \sum_{s \in \mathcal{S}} p_s c(\vx_s) 
    \tag{SND}\label{eq:sndta}\\
    \text{s.t.} \;\; 
    & y_e = 0 \Rightarrow x_{e,s} \leq 0 
    && \forall e \in \mathcal{R},\; s \in \mathcal{S}, \nonumber\\
    & \vx_s \in \mathcal{F}_s 
    && \forall s \in \mathcal{S}, \nonumber\\
    & \vy \in \mathcal{Y}. \nonumber
\end{align}
Here, $\mathcal{R}$ is the set of candidate edges and $\vy$ are the associated design variables. 
The implication constraints link the design variables with the scenario-dependent flows: if a 
candidate edge is not selected, then no flow can be routed through this edge in any scenario. The 
design variables are shared across all scenarios, whereas the flow variables and flow constraints 
are scenario-specific.
}

\revision{
We also consider the penalty-based formulation, where the linking 
constraints are replaced by a nonlinear penalty term:
\begin{align}
    \min_{\vy, \vx} \;\; 
    & \mathbf{r}^\intercal \vy 
    + \sum_{s \in \mathcal{S}} p_s
    \left(
        c(\vx_s)
        + \mu 
        \sum_{z \in \mathcal{Z}} 
        \sum_{e \in \mathcal{R}} 
        \max\{x_{e,s}^z - M^z y_e, 0\}^p
    \right)
    \tag{PB-SND}\label{eq:srndta}\\
    \text{s.t.} \;\; 
    & \vx_s \in \mathcal{F}_s 
    && \forall s \in \mathcal{S}, \nonumber\\
    & \vy \in \mathcal{Y}. \nonumber
\end{align}
This formulation removes the linking constraints from the feasible region. As a result, the LMO can
decompose into independent subproblems for the design variables and the flow variables. The 
flow LMO can be evaluated separately for each scenario, while the design LMO is computed over 
$\mathcal{Y}$. This greatly reduces the cost of the oracle, but weakens the formulation since feasibility with 
respect to the linking constraints is enforced only through the penalty term.
}

\revision{
We use the five network instances from 
\citet{transportnetlibrary}: \textit{Berlin-Friedrichshain}, \textit{Berlin-Tiergarten}, \textit{Berlin-Prenzlauerberg-Center}, 
\textit{BerlinMitte-Center} and \textit{Anaheim}. 
For each network, we generate five instances with different random seeds. In each instance, between $1\%$ and $5\%$ 
of the edges are removed uniformly at random and used as candidate design edges. An instance is 
considered solved if the relative dual gap is below $5\%$. For the penalty-based formulation, we also 
require the maximum violation of the linking constraints to be below $0.01$. The time limit is 
$3600$ seconds per instance.
}

\paragraph{\revision{Simple LMO and MOI LMO.}}

\revision{
We compare the direct formulation \eqref{eq:sndta} with the penalty-based formulation 
\eqref{eq:srndta}. For the direct formulation, the LMO is modeled through MOI and solved with 
SCIP or HiGHS. For the penalty-based formulation, we use either a specialized simple LMO that 
exploits the decomposable structure, or a generic MOI-based LMO. The results 
for $1$ and $5$ scenarios are shown in \cref{subfig:network_1,subfig:network_5}.
}

\revision{
For $1$ scenario, the direct formulation solves more instances within the time limit overall. 
In this regime, the MOI-based LMO for the direct formulation remains small in size, while 
the linking constraints yield substantially stronger relaxations. Specifically, in these experiments, 
the root-node lower bounds are already strong and, during the Frank-Wolfe process, Boscia can 
already identify an integer-feasible incumbent that brings the optimality gap below the tolerance.
}

\revision{
For $5$ scenarios, the MOI-based LMO becomes substantially more expensive because the number of scenario-dependent flow variables increases, 
while the design variables remain coupled across all scenarios through the linking constraints. 
In this setting, the penalty-based formulation with the simple LMO solves the largest number of instances. 
Although the direct formulation yields much stronger relaxations and can certify optimality at the root node, the cost of its LMO becomes prohibitive. 
As a result, Boscia cannot perform enough Frank-Wolfe iterations within the time limit to find an integer-feasible incumbent and close the optimality gap, 
even at the root node. 
In contrast, the cheaper simple LMO allows the penalty-based formulation to find high-quality solutions after processing a moderate number of nodes.
}

\revision{
The penalty-based formulation with MOI-based LMOs performs poorly in both settings. 
Here, the separability of the formulation is not exploited, so the oracle calls remain expensive while the relaxation is much weaker.
This comparison shows that the penalty-based formulation is effective only when paired with a specialized LMO. 
For a small number of scenarios, the tighter direct formulation is preferable; for more scenarios, the cheaper oracle of the penalty-based formulation becomes dominant.
}

\begin{figure}[H]
    \centering
    \subfloat[1 scenario]{\label{subfig:network_1}
        \wenjieplot[width=0.46\textwidth]{experiments/Network/1_scenario_comparison.pdf}
    }\hfill
    \subfloat[5 scenarios]{\label{subfig:network_5}
        \wenjieplot[width=0.46\textwidth]{experiments/Network/5_scenarios_comparison.pdf}
    }
    \caption{\revision{Network design experiments with $1$ and $5$ scenarios. The direct formulation is solved with MOI-based LMOs, while the penalty-based formulation is solved either with a specialized simple LMO or with MOI-based LMOs.}}
    \label{fig:network_scenarios}
\end{figure}

\paragraph{\revision{Frank-Wolfe variants.}}

\revision{
We next compare different Frank-Wolfe variants as node solvers for the penalty-based formulation 
with the specialized simple LMO. The results are shown in \cref{fig:fw_variants}. The comparison 
includes \emph{Standard FW}, \emph{AFW}, \emph{PFW}, \emph{BCG}, and 
\emph{BPCG}. Among these variants, \emph{AFW} solves the largest number of 
instances within the time limit, followed by \emph{BPCG}. One possible explanation is that, 
for the continuous relaxations arising in this problem, away directions appear to provide particularly effective 
descent directions compared with the alternative descent directions used by the other Frank-Wolfe variants.
}

\revision{
Note that \texttt{Standard FW} performs worst in this experiment. This is expected, since standard 
Frank-Wolfe requires more LMO calls overall. Even with the simple LMO, these oracle 
calls remain the dominant cost in the node solves. Variants that can make progress using the current 
active set, such as away-step or pairwise variants, and that can further benefit from lazification, 
reduce the number of LMO calls and perform better in this setting.
}

\begin{figure}[H]
    \centering
    \subfloat[Frank-Wolfe variants]{\label{fig:fw_variants}
        \wenjieplot[width=0.46\textwidth]{experiments/Network/FW_variants.pdf}
    }\hfill
    \subfloat[Line-search methods]{\label{fig:line_searches}
        \wenjieplot[width=0.46\textwidth]{experiments/Network/line_searches.pdf}
    }
    \caption{\revision{Network design experiments for the penalty-based formulation with the specialized simple LMO: (a)~Frank--Wolfe variants as node solvers; (b)~line-search methods.}}
\end{figure}

\paragraph{\revision{Line-search methods.}}

\revision{
    We compare different line-search methods for the penalty-based formulation with the specialized simple LMO.
    The results are shown in \cref{fig:line_searches}.
    \revII{The Adaptive line search performs best, solving the largest fraction of instances within the time limit, followed by BackTracking.
    Both methods replace the global smoothness constant by a local estimate, but differ in the criterion used to accept this estimate.
    The superior performance of Adaptive suggests that, for these instances, its gradient-based acceptance criterion yields more effective local smoothness estimates than the sufficient-decrease condition employed by BackTracking.
    Nevertheless, BackTracking remains competitive and achieves the second-best performance overall.

    The Agnostic strategy, using $\ell=2$ and hence the step-size rule $\gamma_t=2/(t+2)$, achieves the third-best performance despite relying on a predetermined schedule that does not adapt to local information. 
    Its competitive performance can be partly attributed to its low computational overhead.
    The Secant line search performs slightly worse, which may be attributed to the additional derivative evaluations required to locate a stationary point along the search direction.
    ShortStep and GoldenRatio are the weakest strategies in this experiment.
    For ShortStep, we estimate a global Lipschitz constant. This global estimate can be conservative relative to the local smoothness and consequently lead to small step sizes.
    GoldenRatio, in contrast, requires repeated objective evaluations to progressively contract the interval containing a minimizer.}
}

\paragraph{\revision{Branching strategies.}}

\revision{
Figure~\ref{fig:branching_strategies} compares different branching strategies on the network design instances.
The \emph{Most Infeasible} rule performs best overall: it starts solving instances earlier and reaches the largest number of solved instances within the time limit. 
This suggests that, for this benchmark, branching on the most fractional variable offers a good balance between decision quality and computational cost.
The \emph{Largest Most Infeasible Gradient}, \emph{Hierarchy}, and \emph{Pseudocost} branching rules show broadly similar performance. 

The \emph{Largest Most Infeasible Gradient} rule selects a branching variable by considering both its fractional infeasibility and its gradient information.
The \emph{Pseudocost} rule uses information from previous branching decisions to estimate which variable is most likely to improve the bound.
The \emph{Hierarchy} rule selects the branching variable through a sequence of prioritized selection stages. 
In our experiments, we use the default hierarchy: candidates are first filtered by the \emph{Most Infeasible} rule, then refined using \emph{Pseudocost} information, and finally selected according to the \emph{Largest Most Infeasible Gradient} criterion if multiple candidates remain.

Although these strategies use additional information, it does not lead to a clear improvement over the simpler most-infeasible rule. 
One possible reason is that the main difficulty in these instances is driven more by integrality violation than by objective sensitivity or historical branching estimates.
In contrast, \emph{Partial Strong Branching} solves noticeably fewer instances. 
This indicates that the extra time spent evaluating candidate variables may not be worthwhile here, since the stronger local branching decisions do not compensate for the additional overhead. 
Overall, the results suggest that the simple Most Infeasible rule is the most effective and robust choice for this experiment.
}

\begin{figure}[H]
    \centering
    \wenjieplot[width=0.7\textwidth]{experiments/Network/branching_strategies.pdf}
    \caption{\revision{Comparison of different branching strategies.}}
    \label{fig:branching_strategies}
\end{figure}

\subsection{Graph isomorphism problem - creating a self-managed LMO}\label{sec:Examples:GraphIsomorphism}

The Graph Isomorphism Problem (GIP) asks whether two graphs $G_1$ and $G_2$ are structurally 
identical, i.e., whether there exists a permutation of the vertex set that maps the adjacency 
structure of one graph onto that of the other, see \cref{fig:petersen-isomorphic} for a visualization.
If $A$ and $B$ denote their adjacency matrices, 
the graphs are isomorphic if and only if
\[
P A P^\top = B,
\]
for some permutation matrix $P$. An equivalent reformulation expresses the problem as a quadratic 
optimization program over the set of permutation matrices:
\[
\min_{X \in \mathcal{P}(n)} \; \lVert X A - B X \rVert_F^2,
\]
where $\mathcal{P}(n)$ is the set of $n \times n$ permutation matrices and 
$\lVert \cdot \rVert_F$ denotes the Frobenius norm. 

\revision{
    Building on the work of \citet{klus2025continuous}, who propose a Frank-Wolfe-based approximation method, and \citet{xiao2025GI}, who 
    apply Boscia to solve the quadratic formulation exactly, it is worthwhile to revisit the Boscia implementation. 
    In particular, it serves as a good example of a fully self-managed LMO.
}

The continuous problem solved with FW at each node in Boscia is defined as:
\[
\min_{X \in \mathcal{D}(n)\cap [\vl, \vu]} \; \lVert X A - B X \rVert_F^2,
\]
where
\[
\mathcal{D}(n) = \left\{ X \in \mathbb{R}^{n \times n} \;\mid\; X \geq 0,\; \sum\limits_{i=1}^n X_{ij} = 1 \;\forall j \in \{1, \ldots, n\},\; \sum\limits_{j=1}^n X_{ij} = 1 \;\forall i \in \{1, \ldots, n\} \right\} \]
denotes the set of doubly stochastic matrices, i.e., the convex hull of the permutation matrices 
$\mathcal{P}(n)$. The additional box constraints $[\vl,\vu]$ are node-specific bounds imposed by the 
branching process, with $\vl, \vu \in \{0, 1\}^{n\times n}$. 

\begin{figure}[htbp]
\centering
\wenjiefigure{%
% ---------- (a) Classic: original labeling ----------
\begin{tikzpicture}[scale=0.75, every node/.style={transform shape},
    baseline=(current bounding box.center)]
\def\R{2.2}\def\r{1.05}\def\phi{90}

% hidden coordinates = original vertices
\foreach \i in {0,...,4}
\path ({\phi+72*\i}:\R) coordinate (A\i);
\foreach \i in {0,...,4}{
\pgfmathtruncatemacro{\j}{\i+5}
\path ({\phi+72*\i+36}:\r) coordinate (A\j);
}

% edges of the Petersen graph in hidden original coordinates
\foreach \u/\v in {
0/1,1/2,2/3,3/4,4/0,
5/7,7/9,9/6,6/8,8/5,
0/5,1/6,2/7,3/8,4/9}
\draw[thick] (A\u) -- (A\v);

% displayed labels: identity
\foreach \i in {0,...,9}
\node[circle,draw,fill=white,minimum size=10pt,inner sep=1pt] at (A\i) {\i};
\end{tikzpicture}
\hspace{0.8cm}
%
% ---------- (b) Circle layout with permuted labeling ----------
\begin{tikzpicture}[scale=0.75, every node/.style={transform shape},
    baseline=(current bounding box.center)]
\def\R{2.2}\def\phi{90}

% hidden coordinates = original vertices, placed on a circle
\foreach \i in {0,...,9}
\path ({\phi+36*\i}:\R) coordinate (B\i);

% same hidden Petersen edges
\foreach \u/\v in {
0/1,1/2,2/3,3/4,4/0,
5/7,7/9,9/6,6/8,8/5,
0/5,1/6,2/7,3/8,4/9}
\draw[thick] (B\u) -- (B\v);

% displayed labels are permuted:
% old vertex: 0 1 2 3 4 5 6 7 8 9
% new label:  2 7 0 9 5 1 8 4 6 3
\foreach \old/\new in {
0/2, 1/7, 2/0, 3/9, 4/5,
5/1, 6/8, 7/4, 8/6, 9/3}
\node[circle,draw,fill=white,minimum size=10pt,inner sep=1pt] at (B\old) {\new};
\end{tikzpicture}
\hspace{0.8cm}
%
% ---------- (c) Elliptical layout with another permuted labeling ----------
\begin{tikzpicture}[scale=0.9, every node/.style={transform shape},
    baseline=(current bounding box.center)]
\def\R{1.8}
\def\r{0.8}

% hidden coordinates = original vertices
\path (  0:\R) coordinate (C0);
\path (330:\R) coordinate (C5);
\path (300:\R) coordinate (C8);
\path (240:\R) coordinate (C6);
\path (210:\R) coordinate (C9);
\path (180:\R) coordinate (C4);
\path (120:\R) coordinate (C3);
\path ( 60:\R) coordinate (C2);
\path ( 30:\R) coordinate (C1);
\path ( 45:\r) coordinate (C7);

% same hidden Petersen edges
\foreach \u/\v in {
0/1,1/2,2/3,3/4,4/0,
5/7,7/9,9/6,6/8,8/5,
0/5,1/6,2/7,3/8,4/9}
\draw[thick] (C\u) -- (C\v);

% displayed labels are permuted:
% old vertex: 0 1 2 3 4 5 6 7 8 9
% new label:  4 0 8 1 6 9 2 5 3 7
\foreach \old/\new in {
0/4, 1/0, 2/8, 3/1, 4/6,
5/9, 6/2, 7/5, 8/3, 9/7}
\node[circle,draw,fill=white,minimum size=12pt,inner sep=0pt] at (C\old) {\new};
\end{tikzpicture}
}

\caption{\revision{Three isomorphic representations of the Petersen graph adapted from the example in \citet{hasan2017graphettes}. 
All three layouts depict the same abstract graph using different drawings and vertex labelings. 
They preserve the same adjacency structure up to relabeling, as well as the same graph-theoretic invariants, 
illustrating the concept of graph isomorphism.}
\label{fig:petersen-isomorphic}}

\end{figure}

The inclusion of the Graph Isomorphism Problem as a tutorial example serves to 
illustrate several features of \package. First, the problem is naturally 
formulated over the Birkhoff polytope, providing a setting in which convex 
relaxations and combinatorial constraints interact directly. Second, it offers 
a clear example for demonstrating the implementation of a custom 
\lstinline{BirkhoffLMO}, together with the associated oracles and node-bound 
management routines. Finally, the problem highlights the role of user-defined 
callbacks, executed before and after node evaluation, which can significantly 
influence the efficiency of the solution process.

We start with the definition and implementation of the \lstinline{BirkhoffLMO}, 
which encodes the structural information required during the branch-and-bound process.
While we could utilize the \lstinline{ManagedLMO} type in \package, its default settings may not suffice.
As stated earlier, the LMO is stored at tree level and updated for each node. 
In the implementation via the \lstinline{ManagedLMO}, the new bounds are only known to the \lstinline{ManagedLMO} object 
and are \revision{provided} to the underlying LMO explicitly during the \lstinline{compute_extreme_point} call.
In the case of the Birkhoff polytope, this is not computationally efficient because fixing a variable to one eliminates 
a row and a column, requiring reconstruction of the reduced problem for each \lstinline{compute_extreme_point} call.
In such cases, it is necessary to define a customized \lstinline{LMO} that extends the default behavior 
and incorporates structural information such as fixed indices and reduced mappings, ensuring 
that the feasible region is updated in a consistent and structure-aware manner.

The Birkhoff polytope 
admits efficient linear minimization through linear assignment problems, making 
it natural to design a dedicated oracle that explicitly exploits this structure. 
In particular, the \lstinline{BirkhoffLMO} keeps track of indices fixed to one as 
well as the reduced matrix obtained by eliminating the corresponding rows and 
columns.
This representation enables fast oracle computations while preserving 
the node-specific combinatorial structure of the feasible region.

\noindent
\begin{minipage}{\linewidth}
\begin{jllisting}
mutable struct BirkhoffLMO <: FrankWolfe.LinearMinimizationOracle
    dim::Int
    lower_bounds::Vector{Float64}
    upper_bounds::Vector{Float64}
    int_vars::Vector{Int}
    fixed_to_one_rows::Vector{Int}
    fixed_to_one_cols::Vector{Int}
    index_map_rows::Vector{Int}
    index_map_cols::Vector{Int}
    updated_lmo::Bool
    atol::Float64
    rtol::Float64
end
\end{jllisting}
\end{minipage}

The fields 
\lstinline{lower_bounds} and \lstinline{upper_bounds} encode node-specific bounds, 
which at the root node are initialized to zero and one, respectively. For each 
node, the \lstinline{BirkhoffLMO} also records the entries fixed to one. In the 
Birkhoff polytope, such a fixing implies that all other entries in the same row 
and column must be set to zero. To efficiently reflect these reductions, the 
fields \lstinline{index_map_rows} and \lstinline{index_map_cols} store the original 
indices of the remaining rows and columns, allowing the Hungarian algorithm to 
be applied directly to the reduced subproblem without re-indexing the entire 
matrix. Finally, the oracle maintains numerical tolerances (\lstinline{atol}, 
\lstinline{rtol}) and an update flag (\lstinline{updated_lmo}), which facilitate 
robust computations during the optimization process.

After formulating the \lstinline{BirkhoffLMO}, we turn to the implementation of the 
core functionality required to solve the continuous relaxation at each node: the LMO extreme point oracle. 
This is realized by the function \lstinline{compute_extreme_point}. 
\revision{A detailed implementation, together with those of several related functions discussed later in this section, is provided in Appendix~\ref{app:graph_isomorphism}.}

The function takes as input the \lstinline{lmo}, which encodes the structure of the
feasible region, and a direction $D$. It returns an extreme point
minimizing the linear objective induced by~$D$. For the Birkhoff polytope, the extreme points are precisely the permutation
matrices. They can be computed efficiently by solving a linear assignment
problem via the Hungarian algorithm. 
This approach proceeds in three steps.
\begin{enumerate}
    \item We assemble the reduced direction matrix $D_{\text{reduced}}$ by restricting $D$ to
    the active rows and columns (given by \lstinline{index_map_rows} and
    \lstinline{index_map_cols}) and forbidding arcs \revision{whose node-specific upper bounds
    in \lstinline{lmo} are equal to zero.}
    \item We solve the resulting linear assignment problem on $D_{\text{reduced}}$ using
    the Hungarian algorithm to obtain an optimal assignment on the reduced
    direction matrix.
    \item We lift this reduced solution to the full $n\times n$ matrix by
    reinstating the removed rows and columns and restoring any entries fixed to one,
    thereby obtaining the desired extreme point.
\end{enumerate}

When the Decomposition-Invariant Conditional Gradient (DICG) variant is employed 
to solve the node subproblem, two additional oracles are required: one for computing 
an in-face extreme point and another for determining the maximum allowed step size. 
These oracles are essential because, unlike active-set-based Frank-Wolfe methods, DICG does 
not explicitly maintain an active set during its iterations.

\revision{
The function \lstinline{compute_inface_extreme_point} takes as input the feasible region 
\lstinline{lmo}, a search direction $D$, and the current iterate $X$, and returns the 
extreme point satisfying the corresponding in-face constraints. Its implementation mirrors 
the computation of a standard extreme point, with the additional restriction that the update 
remains within the minimal face containing $X$. 
The step-size oracle \lstinline{dicg_maximum_step} takes as input \lstinline{lmo}, the current iterate $X$, and a search 
direction $D$, and computes the largest $\gamma_{\max} \in [0,1]$ such that $X-\gamma D$ remains feasible. It 
iterates over all entries of $X$, computes the maximum step size before each entry reaches the boundary, and returns the minimum of these values.
}

\revision{
Since the \lstinline{BirkhoffLMO} is initialized by the user and dynamically updated by the framework at each node, several interfaces must be implemented to support these updates.
For example, to access the bounds, the \lstinline{get_lower_bound_list}, \lstinline{get_upper_bound_list}, and \lstinline{get_integer_variables} functions must be implemented accordingly. 
To modify the bounds, the \lstinline{set_bound!} function must be implemented.
}

\revision{
In particular, the implementation of the \lstinline{delete_bounds!} function is specifically designed for the Birkhoff polytope. 
The routine constructs the reduced matrix indices based on the \lstinline{fixed_to_one_rows} and \lstinline{fixed_to_one_cols} fields of the LMO. This reduction is performed 
only once at each node. In contrast, under the framework’s default bound management, it is performed at every LMO call, which is computationally inefficient.
}

In addition to the bound-handling routines, safety checks are incorporated to ensure 
the validity of candidate solutions. As an example, the function 
\lstinline{is_linear_feasible} takes the LMO and a vertex as input and verifies that all 
entries respect the bounds currently saved in the LMO, and that both the row and column 
sums are equal to one. Satisfying these conditions certifies that the candidate point 
belongs to the Birkhoff polytope.

Furthermore, additional utility functions are available to help monitor the correctness 
of the solution process. While their use is optional, they can be valuable in practice. 
An example is the function \lstinline{build_LMO_correct}, which verifies that the linear 
minimization oracle has been constructed consistently with the imposed node-specific 
bounds.

In the Graph Isomorphism Problem, if and only if two graphs are isomorphic is the minimum attainable 
objective value exactly zero. Consequently, as soon as a feasible permutation matrix 
is identified that achieves an objective value 
of zero, the branch-and-bound procedure may be terminated immediately without loss of 
correctness. \revision{Moreover, if the lower bound of a node is already strictly positive, none of its child nodes can contain a 
feasible solution certifying an isomorphism. Therefore, no branching is performed.}

The \package{} supports such problem-specific dynamic behavior through user-defined callbacks. In particular, two interfaces are 
provided: \lstinline{branch_callback}, which is executed during branching, and \lstinline{bnb_callback}, which is invoked before the next node is evaluated. 
These callbacks enable users to embed custom logic into the branch-and-bound procedure. For the Graph Isomorphism Problem, this makes early stopping straightforward: 
once an isomorphism is detected, \revision{that is}, a permutation matrix with objective value zero is found, the \lstinline{bnb_callback} can \revision{terminate the solver 
immediately and avoid unnecessary computation. In addition, the \lstinline{branch_callback} can prune nodes whose lower bound is strictly positive.}

Finally, we configure \package{} with customized settings for the Graph Isomorphism 
Problem:  

\noindent
\begin{minipage}{\linewidth}
\begin{jllisting}
lmo = BirkhoffLMO(n, collect(1:(n^2)))
settings = Boscia.create_default_settings()
settings.branch_and_bound[:verbose] = true
settings.branch_and_bound[:bnb_callback] = build_tree_callback()
settings.branch_and_bound[:branch_callback] = build_branch_callback()
settings.frank_wolfe[:variant] = 
            Boscia.DecompositionInvariantConditionalGradient()
settings.frank_wolfe[:line_search] = FrankWolfe.Secant()
settings.frank_wolfe[:lazy] = true
settings.frank_wolfe[:max_fw_iter] = 1000
# Solve the graph isomorphism problem
x, _, result = Boscia.solve(f, grad!, lmo, settings = settings)
\end{jllisting}
\end{minipage}

In this setup, the branch-and-bound procedure is equipped with user-defined 
callbacks that control how nodes are explored and branched, while progress 
statistics are reported periodically. 
On the continuous optimization side, DICG is run with the 
Secant line search, lazy updates, and a cap of $1000$ iterations to balance 
efficiency with solution quality.  

\subsubsection*{\revision{Computational results}}
\revision{
The experiments are based on 12 graph families. Within each family, the number of nodes varies 
from 10 to 500. For each graph, we randomly generate three isomorphic pairs, which serve as the 
test instances. An instance is considered solved if the computed solution satisfies 
$A \approx XBX^\revII{\top}$. The time limit is set to 
$3600$ seconds per instance. Unless stated otherwise, we use 
\texttt{DICG} as the node solver and \texttt{Secant} as the 
line-search method. We also use \texttt{BiasedDepthFirst} as the node selection strategy and leverage the \texttt{star} preprocessing
strategy, which is reported to be most effective for the GI problem on these instances in \citet{xiao2025GI}.
}

\begin{table}[h]
    \centering
    \wenjietable{%
    \footnotesize
    \setlength{\tabcolsep}{2.2pt}%
    \resizebox{\dimexpr\linewidth-12pt\relax}{!}{%
    \begin{tabular}{lrrrrrrrrrrrr}
        \toprule
         & \multicolumn{3}{c}{Self-Managed LMO} & \multicolumn{3}{c}{Simple LMO} & \multicolumn{3}{c}{MOI LMO (SCIP)} & \multicolumn{3}{c}{MOI LMO (HiGHS)} \\
        \cmidrule(lr){2-4}\cmidrule(lr){5-7}\cmidrule(lr){8-10}\cmidrule(lr){11-13}
        Family & \makecell{\%\\solved} & Time & sec/node & \makecell{\%\\solved} & Time & sec/node & \makecell{\%\\solved} & Time & sec/node & \makecell{\%\\solved} & Time & sec/node \\
        \midrule
        latin & \textbf{44\%} & \textbf{207.17} & \textbf{26.19} & 41\% & 232.44 & 30.83 & 32\% & 603.76 & 102.50 & 32\% & 757.34 & 102.90 \\
        Lattice & \textbf{57\%} & \textbf{321.82} & \textbf{22.68} & \textbf{57\%} & 394.86 & 26.53 & 29\% & 1148.32 & 137.23 & 43\% & 1708.34 & 129.62 \\
        paley\_power & \textbf{62\%} & \textbf{120.72} & \textbf{16.90} & 58\% & 171.27 & 26.39 & 42\% & 480.56 & 91.05 & 50\% & 653.97 & 89.94 \\
        paley\_prime & \textbf{71\%} & \textbf{96.62} & \textbf{20.21} & \textbf{71\%} & 130.50 & 24.78 & \textbf{71\%} & 539.57 & 93.52 & \textbf{71\%} & 525.96 & 87.58 \\
        Triangular & \textbf{57\%} & \textbf{450.58} & \textbf{27.67} & 52\% & 548.05 & 37.66 & 29\% & 1707.85 & 175.24 & 43\% & 1955.47 & 154.36 \\
        CHH\_cc & \textbf{33\%} & \textbf{925.94} & \textbf{45.20} & 29\% & 1034.62 & 95.57 & 21\% & 2291.04 & 409.97 & 21\% & 2517.13 & 401.28 \\
        tnn & \textbf{31\%} & \textbf{961.62} & \textbf{24.84} & 26\% & 1231.85 & 45.19 & 17\% & 2225.80 & 268.16 & 14\% & 2816.68 & 319.60 \\
        cfi & \textbf{0\%} & \textbf{3600.00} & \textbf{268.51} & \textbf{0\%} & \textbf{3600.00} & 749.44 & \textbf{0\%} & \textbf{3600.00} & 1441.19 & \textbf{0\%} & \textbf{3600.00} & 1200.00 \\
        exact & \textbf{100\%} & \textbf{4.67} & \textbf{2.74} & 95\% & 7.05 & 4.09 & 86\% & 88.73 & 39.13 & 86\% & 133.56 & 60.62 \\
        iso\_r01N & \textbf{86\%} & \textbf{2.28} & \textbf{1.87} & \textbf{86\%} & 2.31 & 1.87 & \textbf{86\%} & 9.06 & 6.12 & \textbf{86\%} & 18.60 & 12.70 \\
        sts & \textbf{57\%} & \textbf{343.10} & \textbf{65.98} & 52\% & 616.47 & 118.17 & 0\% & 3600.00 & 533.18 & 0\% & 3600.00 & 383.14 \\
        usr & \textbf{69\%} & \textbf{79.93} & \textbf{17.34} & \textbf{69\%} & 95.46 & 20.26 & 17\% & 2452.95 & 210.22 & 6\% & 3283.23 & 256.71 \\
        Summary & \textbf{55\%} & \textbf{179.14} & \textbf{21.31} & 52\% & 220.02 & 30.80 & 34\% & 817.18 & 142.37 & 35\% & 1037.44 & 153.99 \\
        \bottomrule
        \end{tabular}
    }%
    }
    \caption{\revision{Performance summary by graph family. Reported times are the geometric mean of per-instance total solving times,
    after shifting each value by one second. Timed-out instances are included. The best result in each row appears in bold.}}
    \label{tab:gi_results_summary}
\end{table}

\revision{
In this experiment, we compare three LMO implementations: the self-managed LMO,
the simple LMO, and the MOI-based LMO using either SCIP or HiGHS as the MIP solver. The results are summarized 
in \cref{tab:gi_results_summary}. Overall, the self-managed LMO performs best: it solves the largest number 
of instances and attains the lowest shifted geometric mean running time across all graph families, except for 
the \texttt{cfi} family, for which no LMO implementation solves any instance within the time limit. In contrast, 
the MOI-based LMOs perform poorly on most families. This is because the Hungarian algorithm is much more efficient than general MIP methods, and 
when using DICG as the Frank-Wolfe variant, two LMO calls are required at each iteration, which further increases the computational overhead. 
Among the two MOI-based variants, SCIP generally outperforms HiGHS as a backend on most graph families.
}

\revision{
For the \texttt{sts} and \texttt{usr} families, the performance gap between the Hungarian-algorithm-based LMOs (the self-managed LMO and the simple LMO)
and the MOI-based LMOs is substantially larger than for the other families. This can be explained by the fact
that, for most instances in these two families, more than half of the variables are fixed to zero by the \texttt{star} preprocessing
step. Although this significantly reduces the size of the feasible region, the linear objectives solved within the LMO
often admit a much larger number of optimal solutions. As a result, the Hungarian-algorithm-based LMOs and the MOI-based LMO may return different
optimal extreme points and therefore follow different solution paths. For these two families, the extreme points returned
by the Hungarian algorithm appear to be more favorable for progressing toward an integer solution.
}

\begin{figure}[H]
    \centering
    \subfloat[\texttt{sts}]{\label{fig:GI_sts}
        \wenjieplot[width=0.46\textwidth]{experiments/GI/solved_vs_time_sts.pdf}
    }\hfill
    \subfloat[\texttt{usr}]{\label{fig:GI_usr}
        \wenjieplot[width=0.46\textwidth]{experiments/GI/solved_vs_time_usr.pdf}
    }
    \caption{\revision{Graph isomorphism experiments: cumulative number of solved instances over time for different LMO implementations on the \texttt{sts} and \texttt{usr} families.}}
    \label{fig:GI_sts_usr_pair}
\end{figure}

\revision{In \cref{tab:gi_results_lmo_time}, we report more detailed timing results for the Hungarian-algorithm-based LMOs. 
The self-managed LMO consistently outperforms the simple LMO, with the advantage becoming particularly pronounced on large 
or difficult instance families, such as \texttt{cfi} and \texttt{tnn}. 
This behavior is expected: on hard instances, either a larger number of nodes must be explored or the underlying graphs are larger, 
so repeatedly reconstructing the subproblem at each iteration incurs a substantial overhead for the simple LMO. 
In contrast, the self-managed LMO stores the reduced indices and reuses this information across calls, thereby 
significantly reducing the per-node computational cost.}

\begin{table}[h]
    \centering
    \wenjietable{%
    \footnotesize
    \setlength{\tabcolsep}{2.2pt}%
    % Few columns: full-\linewidth resizebox blows up font height; cap width (same style as first table, smaller target)
    \resizebox{0.48\linewidth}{!}{%
    \begin{tabular}{lrrrr}
        \toprule
         & \multicolumn{2}{c}{Self-Managed LMO} & \multicolumn{2}{c}{Simple LMO} \\
        \cmidrule(lr){2-3}\cmidrule(lr){4-5}
        Family & \makecell{\%\\solved} & ms/call & \makecell{\%\\solved} & ms/call \\
        \midrule
        latin & \textbf{44\%} & \textbf{22.156} & 41\% & 150.352 \\
        Lattice & \textbf{57\%} & \textbf{12.254} & \textbf{57\%} & 50.195 \\
        paley\_power & \textbf{62\%} & \textbf{6.237} & 58\% & 58.418 \\
        paley\_prime & \textbf{71\%} & \textbf{5.910} & \textbf{71\%} & 33.648 \\
        Triangular & \textbf{57\%} & \textbf{14.202} & 52\% & 63.198 \\
        CHH\_cc & \textbf{33\%} & \textbf{31.199} & 29\% & 176.619 \\
        tnn & \textbf{31\%} & \textbf{9.707} & 26\% & 81.622 \\
        cfi & \textbf{0\%} & \textbf{178.307} & \textbf{0\%} & 4571.771 \\
        exact & \textbf{100\%} & \textbf{0.609} & 95\% & 43.263 \\
        iso\_r01N & \textbf{86\%} & \textbf{0.172} & \textbf{86\%} & 10.545 \\
        sts & \textbf{57\%} & \textbf{30.648} & 52\% & 674.772 \\
        usr & \textbf{69\%} & \textbf{1.948} & \textbf{69\%} & 55.521 \\
        Summary & \textbf{55\%} & \textbf{9.667} & 52\% & 93.959 \\
        \bottomrule
        \end{tabular}
    }%
    }

    \caption{\revision{Performance summary by graph family for the Hungarian-algorithm-based LMOs.
    Reported LMO times (\texttt{ms/call}) are geometric means shifted by one millisecond.
    Timed-out instances are included. The best result in each row appears in bold.}}
    \label{tab:gi_results_lmo_time}
\end{table}

\subsection{Optimal design of experiments - restricted function domains}\label{sec:Examples:OEDP}

The \emph{Optimal Experiment Design Problem (OEDP)} presents a compelling case study for mixed-integer convex optimization, 
as it naturally combines the discrete nature of experiment selection with convex information measures.
In OEDP, we are given a matrix $A \in \mathbb{R}^{m \times n}$ where each row represents a potential experiment, 
and our goal is to select a subset of $N$ experiments that maximizes some measure of information about the parameter space.
Note that the number of experiments $m$ is typically much larger than the number of experiments to select $N$.
Furthermore, we assume that the matrix $A$ is full-rank, i.e., $\rank(A) = n$, and $n \leq N$.
This leads to a pure-integer nonlinear program where the integer variables indicate which experiments to include,
while the objective function involves various information measures that are convex in the continuous relaxation.

Often, the ultimate aim is to fit parameters \revision{of} a linear model given by \revision{$A$} and match the model to the experiment data. 
In this setting, maximizing information is the same as minimizing the variance of the parameter estimates, 
as the inverse of the \emph{information matrix} $X(\vx) = A^T \text{diag}(\vx) A$ is the so-called \emph{dispersion matrix} 
encoding the covariances of the parameters.
Thus, the choice of information measure determines what aspect of parameter estimation we want to optimize.
The A-criterion, defined as $\Tr(X(\vx)^{-1})$, minimizes the trace of the dispersion matrix,
effectively minimizing the average variance of the parameter estimates.
The D-criterion, formulated as $-\log\det(X(\vx))$, maximizes the determinant of the information matrix,
which corresponds to minimizing the volume of the confidence ellipsoid around the parameter estimates.
Other important criteria include the E-criterion $\lambda_{\min}(X(\vx))$, which maximizes the minimum eigenvalue of the information matrix, 
thereby minimizing the worst possible variance \revision{for any normalized linear combination of parameters}.
Many other information measures exist; here, we focus on the A- and D-criteria.

\begin{minipage}[t]{0.48\textwidth}
\textbf{A-criterion OEDP:}
\begin{align*}
\min_{\vx} \quad & \Tr\left((A^T \text{diag}(\vx) A)^{-1}\right) \\
\text{s.t.} \quad & \sum_{i=1}^m x_i = N \\
& 0 \leq x_i \leq u_i, \quad i = 1, \ldots, m \\
& \vx \in \mathbb{Z}^m 
\end{align*}
\end{minipage}
\hfill
\begin{minipage}[t]{0.48\textwidth}
\textbf{D-criterion OEDP:}
\begin{align*}
\min_{\vx} \quad & -\log\det\left(A^T \text{diag}(\vx) A\right) \\
\text{s.t.} \quad & \sum_{i=1}^m x_i = N \\
& 0 \leq x_i \leq u_i, \quad i = 1, \ldots, m \\
& \vx \in \mathbb{Z}^m 
\end{align*}
\end{minipage}

\noindent
In both formulations, $\vx = (x_1, \ldots, x_m)^T$ represents the number of times each experiment is selected, 
$A \in \mathbb{R}^{m \times n}$ is the design matrix where each row corresponds to a potential experiment,
$N$ is the total budget (number of experiments to select) and 
$u_i$ are upper bounds on the number of times experiment $i$ can be selected.

\revision{Observe that both objective functions are only well defined for positive definite matrices.
Thus, }a common approach for OEDP is based on conic programming formulations, see \citet{sagnol2011computing,sagnol2015computing,coey2022solving}.
The involved cones, however, cannot be represented easily in most conic solvers.
A notable exception is the \texttt{Hypatia.jl} \citep{CoeyLubinVielma2020} solver which implements many natural conic formulations.
\revision{Both problems can be reformulated using semi-definite programming constraints. Interior-point methods like Mosek 
\citep{aps2022mosek} can then be used within branch-and-bound schemes, for example using Mosek as conic solver within SCIPSDP \citep{SCIPSDPGithub,gally2018framework}.}
\revII{Other interior point solvers such as Knitro \citep{byrd2006knitro} and Ipopt \citep{wachter2006implementation} can solve continuous 
relaxations of OEDP under the D-criterion \citep{ponte2023branch,pillai2025computing}.}
An \revision{another} alternative consists in reformulating the nonlinearities with SOCP \citep{sagnol2015computing} constraints, this approach however struggles to scale, 
see \citet{2023_HendrychBesanconPokutta_Optimalexperimentdesign} for a detailed discussion and extensive computational results.

Tackling the above formulations with first-order methods is challenging as both objective functions, and their gradients, 
are only well-defined if the information matrix $X(\vx)$ is positive definite.
However, $X(\vx)$ is not positive definite for all $\vx$ in the feasible region.
Therefore, we cannot start Frank-Wolfe at an arbitrary point of the feasible region.

On the topic of feasible region, notice that it is simply a scaled and truncated probability simplex.
\begin{align}
    S = \left\{\vx \in \mathbb{R}^m : 0 \leq \vx \leq \vu, \sum_{i=1}^m x_i = N\right\}
    \end{align}
where $N$ is the budget for the experiments, $\vu$ are upper bounds, and $m$ is the number of candidate experiments.
That is, the associated LMO is computationally inexpensive and can be solved as a continuous knapsack by sorting gradient entries.

A critical aspect of solving OEDP with the Boscia framework is the proper definition of the domain oracle 
and the ability to generate valid starting points for the subproblems at node level.
The domain oracle determines whether a given point lies within the domain of the objective function, 
i.e., if the activated experiments provide sufficient information for parameter estimation.
From a linear algebra perspective, this is equivalent to the information matrix $X(\vx)$ being positive definite. 
\revision{A detailed implementation is provided in Appendix~\ref{app:optimal_design_problem}, together with those of the related functions discussed later in this section.}

Note that the domain oracle should also be supplied to the line search chosen.
Both the Secant and the Adaptive line search in \texttt{FrankWolfe.jl} can receive domain oracles and 
compute a step size with respect to the domain.

After branching, the active set of a new node might not define a domain-feasible point.
In this case, a domain-feasible point must be provided that respects the new bound constraints,
if possible. 
Otherwise, we assume that the new node is infeasible and it will be pruned.

\revision{First, the domain point routine} computes a set of $n$ linearly independent rows of the experiment matrix $A$ respecting the current upper
bounds. 
Next, it iteratively adds the experiments to the point while respecting the budget constraint $\sum x_i = N$ 
and the upper bounds.
If no domain-feasible point can be computed, \lstinline{nothing} is returned and the node is pruned.

Given a domain-feasible point, we can generate a new starting point by solving a projection problem with 
Frank-Wolfe.
This projection minimizes the distance to the domain-feasible point \revision{using the current point after branching as initial point in the FW. The aim is to produce a good warm start point that is domain-feasible.}
Observe that we do not aim for optimality \revision{of the projection problem} but only want the solution to be in the domain of the original objective (and not on its boundary).

The callback ensures that the algorithm continues for several iterations after finding a domain-feasible point,
so that we do not start with a point on the boundary of the domain.
This can otherwise cause issues in the further solving process.

Finally, we configure \package{} with the appropriate settings and solve the optimization problem:

\noindent
\begin{minipage}{\linewidth}
\begin{jllisting}[language=julia, style=jlcodestyle]{}
# Configure Boscia settings
settings = Boscia.create_default_settings()
settings.branch_and_bound[:verbose] = true
settings.domain[:active_set] = copy(active_set)
settings.domain[:domain_oracle] = domain_oracle
settings.domain[:find_domain_point] = domain_point
settings.heuristic[:hyperplane_aware_rounding_prob] = 0.7
settings.frank_wolfe[:line_search] = FrankWolfe.Secant(
    domain_oracle=domain_oracle
)
settings.frank_wolfe[:lazy] = true

# Solve the A-optimal design problem
x_a, _, _ = Boscia.solve(f_a, grad_a!, lmo, settings=settings)

# Solve the D-optimal design problem  
x_d, _, _ = Boscia.solve(f_d, grad_d!, lmo, settings=settings)
\end{jllisting}
\end{minipage}

\subsubsection*{\revision{Computational results}}

\revision{
    
In the following, we report computational results for several configurations of the Boscia framework applied to the OEDP. In particular, we investigate the effects of different lazification strategies, Frank-Wolfe tolerances, and primal heuristic strategies.
We follow the experimental setup of \citet{2023_HendrychBesanconPokutta_Optimalexperimentdesign}. We consider both the A- and D-optimality criteria, with the number of candidate experiments given by
$m \in \{50,60,80,100,120\}.$
For each value of \(m\), we consider \(n = \left\lfloor \frac{m}{10} \right\rfloor\) and \(n = \left\lfloor \frac{m}{4} \right\rfloor\) parameters. The experimental budget is set to \(N = \lfloor 1.5n \rfloor \). 
\revII{For each resulting parameter configuration ($m$, $n$, $N$), we generate five random instances.}
All lower bounds are fixed at zero, whereas the upper bounds are sampled uniformly at random from the interval \([1,m/10]\).
We consider both independent and correlated experimental data\revision{, both generated randomly}.
Throughout the experiments, we use the default \texttt{Secant} line-search method \citep{hendrych2025secant}, which \citet{2023_HendrychBesanconPokutta_Optimalexperimentdesign} report as the most robust line-search strategy for OEDP instances.
In this experimental setup, an instance is considered solved once the absolute gap falls below \(10^{-6}\).

\paragraph{\revision{Lazification and shadow set.}}

In \cref{fig:OED_lazy}, we compare the performance of different lazification strategies for BPCG. Specifically, we consider 
lazification of the LMO calls both with and without the use of the \revision{shadow set}. The shadow set explicitly stores 
all extreme points generated by the LMO across the nodes of the branch-and-bound tree. As shown in \cref{fig:OED_lazy}, 
lazified BPCG with the shadow set achieves the best overall performance. This behavior is expected, as reusing 
previously generated extreme points can substantially reduce the number of LMO calls. However, explicitly maintaining 
and searching the shadow set introduces additional storage and lookup overhead, which may become significant as the 
branch-and-bound tree grows. Consequently, for some instances, this overhead is not fully offset by the reduction in LMO 
calls, and lazification without the shadow set can yield better performance. This trade-off may also explain why, for 
the independent-design instances under the D-optimality criterion, lazified BPCG with the shadow set solves slightly 
fewer instances than the other two configurations.

Moreover, we observe that, for OEDP, DICG still outperforms lazified BPCG with the shadow set in terms of the number of solved instances. 
\revII{Although DICG does not solve all instances to optimality within the one-hour time limit, the mean absolute gaps reported in
Table~\ref{tab:abs_gap_by_criterion_1e-2_optimal_by_criterion} remain relatively small.
Even for the more challenging configurations with $n=\lfloor m/4\rfloor$, the mean absolute gaps remain below $2.03\times10^{-4}$ and $1.27\times10^{-2}$ 
for A- and D-optimal design with correlated data, respectively, and below $2.46\times10^{-2}$ for both criteria with independent data.
These results suggest that Boscia can still identify relatively high-quality solutions even when optimality cannot be certified within the time limit, as also observed by
\citet{2023_HendrychBesanconPokutta_Optimalexperimentdesign}.}
This can be explained by the fact that DICG achieves more progress per iteration than BPCG and therefore requires fewer iterations 
to close the optimality gap at each node. 
Although each DICG iteration requires two additional oracle calls, the resulting progress appears to compensate for the additional computational cost.
}

\begin{figure}[H]
    \centering
    \wenjieplot[width=0.7\textwidth]{experiments/OED/lazy.pdf}
    \caption{\revision{Lazification and shadow set for OED.}}
    \label{fig:OED_lazy}
\end{figure}

\paragraph{\revision{Different FW epsilons.}}

\revision{
    In Boscia, the parameter \texttt{FW\_epsilon} controls the Frank-Wolfe gap tolerance at the root node. By default, 
    this tolerance decreases as the branch-and-bound tree grows according to a prescribed decay factor. A tighter value 
    of \texttt{FW\_epsilon} generally yields a stronger lower bound at the root node. However, it also increases the 
    computational effort required at each node and may therefore slow down the exploration of the tree. Conversely, a 
    looser value of \texttt{FW\_epsilon} allows nodes to be processed more rapidly, but may produce weaker lower bounds 
    during the early stages of the search. The resulting trade-off is generally problem-dependent.

    The effect of different Frank-Wolfe tolerances on the OEDP instances is shown in \cref{fig:OED_epsilons}\revII{; detailed numerical results are given in Appendix~\ref{app:oed_fw_tolerances}}. 
    \revII{
    For both D-optimal settings, the default value
    $\texttt{FW\_epsilon}=10^{-2}$ performs best overall in terms of both the number of solved instances and the mean absolute gap across the grouped instances.
    For the A-optimal independent instances, $\texttt{FW\_epsilon}=5\times10^{-3}$ solves the largest number of instances, while the default value $10^{-2}$ yields a slightly smaller mean absolute gap overall.
    In contrast, for the A-optimal correlated instances, the choice of tolerance has a much more pronounced effect: tighter values of \texttt{FW\_epsilon} consistently lead to better overall
    performance in terms of the number of solved instances and the mean absolute gap.}
    
    This behavior can be attributed to the primal heuristic employed in Boscia. A tighter tolerance requires more 
    Frank-Wolfe iterations before the stopping criterion is satisfied and therefore generates more vertices through 
    LMO calls. These vertices are subsequently used by the primal heuristic to obtain valid upper bounds. As shown in 
    \cref{tab:best_incumbent_a_correlated}, this leads to substantially earlier identification of \revII{high-quality} integer 
    solutions. Such solutions provide the best possible upper bound and enable more effective pruning during the early 
    stages of the branch-and-bound search, preventing the 
    tree from growing excessively large and substantially improving overall performance in the 
    A-optimal correlated setting.
    \revII{}
    }

\begin{figure}[H]
    \centering
    \wenjieplot[width=0.7\textwidth]{experiments/OED/FW_epsilons.pdf}
    \caption{\revision{Different FW epsilons for OED.}}
    \label{fig:OED_epsilons}
\end{figure}

\begin{table}[htbp]
    \centering
    \wenjietable{%
    \footnotesize
    \setlength{\tabcolsep}{4pt}%
    \resizebox{\dimexpr\linewidth-12pt\relax}{!}{%
    \begin{tabular}{rrrrrrrrrrrrrr}
        \toprule
        \# exp. & \# par. & \multicolumn{2}{c}{5e-2} & \multicolumn{2}{c}{1e-2} & \multicolumn{2}{c}{5e-3} & \multicolumn{2}{c}{1e-3} & \multicolumn{2}{c}{5e-4} & \multicolumn{2}{c}{1e-4} \\
        \cmidrule(lr){3-4}\cmidrule(lr){5-6}\cmidrule(lr){7-8}\cmidrule(lr){9-10}\cmidrule(lr){11-12}\cmidrule(lr){13-14}
         &  & \makecell{\%\\solved} & Time & \makecell{\%\\solved} & Time & \makecell{\%\\solved} & Time & \makecell{\%\\solved} & Time & \makecell{\%\\solved} & Time & \makecell{\%\\solved} & Time \\
         \midrule
         50 & 5 & \textbf{100\%} & 8.87 & \textbf{100\%} & 4.80 & \textbf{100\%} & 5.02 & \textbf{100\%} & 4.31 & \textbf{100\%} & \textbf{3.96} & \textbf{100\%} & 4.19 \\
         50 & 12 & 0\% & 1772.95 & \textbf{80\%} & 119.99 & \textbf{80\%} & 133.14 & \textbf{80\%} & 45.12 & \textbf{80\%} & 47.55 & \textbf{80\%} & \textbf{34.24} \\
         60 & 6 & \textbf{100\%} & 251.26 & \textbf{100\%} & 19.22 & \textbf{100\%} & 8.97 & \textbf{100\%} & 5.66 & \textbf{100\%} & 4.82 & \textbf{100\%} & \textbf{4.55} \\
         60 & 15 & 0\% & 2218.49 & 0\% & 745.33 & 0\% & 270.36 & \textbf{40\%} & 817.11 & \textbf{40\%} & \textbf{236.28} & 20\% & 1193.06 \\
         80 & 8 & 0\% & 2496.78 & \textbf{100\%} & 332.31 & \textbf{100\%} & 87.71 & \textbf{100\%} & 10.52 & \textbf{100\%} & 8.89 & \textbf{100\%} & \textbf{5.48} \\
         80 & 20 & 0\% & 2167.44 & 0\% & 1744.10 & 0\% & \textbf{1594.35} & 0\% & 1765.74 & 0\% & 1604.75 & 0\% & 2003.63 \\
         100 & 10 & 0\% & 2140.65 & 0\% & 2482.40 & 0\% & 1978.29 & 20\% & 462.22 & 60\% & 165.99 & \textbf{100\%} & \textbf{93.81} \\
         100 & 25 & 0\% & 2717.18 & 0\% & 2631.04 & 0\% & 1535.33 & 0\% & 2103.16 & 0\% & \textbf{1468.00} & 0\% & 1704.56 \\
         120 & 12 & 0\% & 1894.94 & 0\% & 2864.99 & 0\% & 2362.60 & 0\% & 1780.70 & 0\% & 895.56 & 0\% & \textbf{123.51} \\
         120 & 30 & 0\% & 2532.49 & 0\% & \textbf{2015.56} & 0\% & 2313.45 & 0\% & 2084.07 & 0\% & 2320.49 & 0\% & 3144.61 \\
         \bottomrule
        \end{tabular}
    }%
    }
    \caption{\revision{Time is the mean time (s) to reach the best incumbent (after which the incumbent no longer improves) for A-optimal
    design with correlated data. The best result in each row is shown in bold.}}
    \label{tab:best_incumbent_a_correlated}
\end{table}

\begin{table}[htbp]
    \centering
    \resizebox{0.85\linewidth}{!}{%
    \begin{tabular}{ccrrrrrrrr}
    \toprule
     &  & \multicolumn{4}{c}{Correlated} & \multicolumn{4}{c}{Independent} \\
    \cmidrule(lr){3-6}\cmidrule(lr){7-10}
    \# exp. & \# par. & \multicolumn{2}{c}{A-optimal} & \multicolumn{2}{c}{D-optimal} & \multicolumn{2}{c}{A-optimal} & \multicolumn{2}{c}{D-optimal} \\
    \cmidrule(lr){3-4}\cmidrule(lr){5-6}\cmidrule(lr){7-8}\cmidrule(lr){9-10}
     &  & \makecell{\%\\solved} & Abs.\ gap & \makecell{\%\\solved} & Abs.\ gap & \makecell{\%\\solved} & Abs.\ gap & \makecell{\%\\solved} & Abs.\ gap \\
     \midrule
     50 & 5 & \textbf{100\%} & 8.39e-07 & \textbf{100\%} & 3.53e-05 & \textbf{100\%} & 2.65e-05 & \textbf{100\%} & \textbf{0.00e+00} \\
     50 & 12 & 80\% & 9.69e-06 & \textbf{100\%} & \textbf{1.07e-07} & 60\% & 3.51e-04 & \textbf{100\%} & 5.90e-07 \\
     60 & 6 & \textbf{100\%} & 9.00e-07 & \textbf{100\%} & 5.16e-08 & \textbf{100\%} & 2.41e-06 & \textbf{100\%} & \textbf{4.84e-09} \\
     60 & 15 & 20\% & 6.57e-05 & \textbf{100\%} & \textbf{7.87e-07} & 0\% & 5.98e-03 & 40\% & 1.94e-03 \\
     80 & 8 & \textbf{100\%} & 9.76e-07 & \textbf{100\%} & 1.94e-07 & \textbf{100\%} & 4.99e-07 & \textbf{100\%} & \textbf{1.41e-07} \\
     80 & 20 & 0\% & \textbf{1.42e-04} & 0\% & 5.72e-03 & 0\% & 1.20e-02 & 0\% & 7.34e-03 \\
     100 & 10 & \textbf{100\%} & 9.99e-07 & \textbf{100\%} & \textbf{4.78e-07} & 60\% & 3.24e-04 & \textbf{100\%} & 9.33e-07 \\
     100 & 25 & 0\% & \textbf{2.03e-04} & 0\% & 9.86e-03 & 0\% & 1.74e-02 & 0\% & 1.36e-02 \\
     120 & 12 & 0\% & 7.59e-06 & \textbf{80\%} & \textbf{3.73e-06} & 0\% & 2.03e-03 & 40\% & 4.68e-04 \\
     120 & 30 & 0\% & \textbf{1.57e-04} & 0\% & 1.27e-02 & 0\% & 2.46e-02 & 0\% & 1.75e-02 \\
     \bottomrule
    \end{tabular}%
    }
    \caption{\protect\revII{Abs.\ gap is the mean absolute gap for A- and D-optimal design with correlated and independent data.
    The best result in each row is shown in bold.}}
    \label{tab:abs_gap_by_criterion_1e-2_optimal_by_criterion}
\end{table}

\paragraph{\revision{Heuristic strategies.}}

\revision{
    Boscia provides several built-in primal heuristics and also supports user-defined heuristics. In our OED experiments, we compare three built-in heuristics, two custom heuristics, and a baseline without 
    heuristics.
    The \texttt{Rounding} heuristic rounds the \revII{non-integer} components of the relaxation solution to the nearest integer. 
    The \texttt{Hyperplane-Aware Rounding} heuristic additionally repairs the rounded point so that the experimental-budget constraint is satisfied. 
    The \texttt{Follow-Gradient} heuristic first follows the objective gradient at the current relaxation solution and calls the LMO to obtain an extreme point.
    It then follows the objective gradient at each newly generated extreme point until a point is repeated or the prescribed number of iterations is reached.
    The two custom heuristics are the \texttt{Fedorov} and the \texttt{Simple-Randomized-Rounding} heuristics, adapted from \citet{fedorov1972theory} and \citet{lamperski2025simple}. 
    The \texttt{Fedorov} heuristic starts from the current incumbent and computes a removal score for each selected experiment. 
    Starting from the lowest-scoring removal candidate, it transfers one unit of mass to other design points with remaining capacity and stops after finding the first improving integer solution.
    The \texttt{Simple Randomized Rounding} heuristic rounds each relaxed component \(x_i\) to \(\lceil x_i\rceil\) with probability \(x_i-\lfloor x_i\rfloor\), and to \(\lfloor x_i\rfloor\) otherwise. 
    Sampling is repeated until the experimental-budget constraint is satisfied or the prescribed number of trials is reached.

    As shown in \cref{fig:OED_heuristics}, the three rounding-based heuristics generally solve more instances than the baseline, with the largest improvement for independent D-optimal design. 
    These heuristics exploit the fractional solution directly and can therefore generate competitive integer incumbents at low computational cost. The stronger effect for independent D-optimal design 
    suggests that its relaxation solutions provide particularly informative selection weights.

    In contrast, \texttt{Follow-Gradient} searches for improved incumbents among extreme points generated by following the gradient direction. When the relaxation solution lies in the interior of the feasible 
    region, these extreme points may be far from it and may therefore discard useful allocation information contained in the fractional design.
    The \texttt{Fedorov} heuristic can even reduce overall performance: it evaluates a sequence of local one-unit exchanges but accepts only the first improvement, so the resulting incumbent gain may 
    be insufficient to offset the computational overhead of the exchange search.

    Moreover, \cref{fig:OED_heuristics_80} shows the evolution of the absolute gap over time for instances with 80 variables; similar behavior is observed for most solved instances. 
    For most of these instances, the rounding-based heuristics identify a high-quality incumbent already at the root node, whereas \texttt{Fedorov} and \texttt{Follow-Gradient} either fail to find a good root-node incumbent or produce a weaker one. 
    This suggests that the relaxed solution obtained at the root node is already close to a good integer solution, making rounding-based strategies effective. 
    In contrast, \texttt{Follow-Gradient} searches among extreme points, which may move farther away from the relaxed solution, while \texttt{Fedorov} performs 
    only local one-unit exchanges and is therefore less effective here. After the early phase, the baseline and most rounding-based heuristics eventually
    identify optimal integer solutions within similar time ranges, which explains why their overall performance profiles remain close. 
    Due to the higher cost of the \texttt{Fedorov} heuristic, however, it typically takes longer to close the optimality gap.
    }   

\begin{figure}[H]
    \centering
    \wenjieplot[width=0.7\textwidth]{experiments/OED/heuristics.pdf}
    \caption{\revision{Comparison of different heuristics for OEDP.}}
    \label{fig:OED_heuristics}
\end{figure}

\begin{figure}[H]
    \centering
    \wenjieplot[width=0.7\textwidth]{experiments/OED/compare_80_8_mean.pdf}
    \caption{\revII{Comparison of different heuristics for OEDP instances with \((m,n,N)=(80,8,12)\). Each curve shows the average absolute gap over five instances in the corresponding setting.}}
    \label{fig:OED_heuristics_80}
\end{figure}

\subsection{\revision{Practical guidelines from the ablation studies}}
\label{sec:Examples:Guidelines}

\revision{
The ablation studies in the three examples above suggest several practical guidelines for using \package{} effectively. 
First, the choice of LMO is often the most important implementation decision. As observed in the network design and graph isomorphism experiments, a customized LMO that exploits an 
efficient algorithm for the underlying linear optimization problem is generally preferable to a generic MIP-solver-based LMO. 
Moreover, if LMO requires internal structure for an efficient extreme point computation, as in the graph isomorphism problem, then a self-managed LMO can be the most effective choice. 
This requires implementing additional interfaces for bound handling and node-level updates, but it can substantially reduce the cost of repeated oracle calls.

Second, the choice of Frank-Wolfe variant should balance per-iteration progress against the cost of oracle calls. If the LMO is relatively cheap and the node subproblems can often be solved 
to sufficient accuracy within the prescribed maximum number of Frank-Wolfe iterations, then DICG can be advantageous, since it typically provides strong progress at each iteration. 
However, DICG requires two additional oracles at each iteration. When oracle calls are expensive, this additional cost may worsen the overall wall-clock performance. 
In such cases, lazified active-set variants, possibly combined with a shadow set, can be more favorable because they can reuse previously generated vertices and reduce the number of 
expensive oracle calls.

Third, the Frank-Wolfe tolerance should be chosen carefully. A tighter initial value of \texttt{FW\_epsilon} can lead to stronger lower bounds and may also generate more useful vertices for 
primal heuristics. This can help identify high-quality incumbents early and prune more nodes in the branch-and-bound tree. 
However, a tighter tolerance also requires more Frank-Wolfe iterations at each node, which can increase node-processing time and lead to worse overall performance if the additional accuracy 
is not offset by stronger pruning.

Finally, the choice of primal heuristic should take both solution quality and computational overhead into account. The default \texttt{Rounding} heuristic provides a good balance 
between cost and effectiveness, especially when the relaxed solution is already close to a high-quality integer solution. More expensive custom heuristics can be useful, but they should not 
necessarily be called after every node evaluation. For example, in the optimal experiment design experiments, the \texttt{Fedorov} heuristic is relatively costly; therefore, using it with 
probability one after every node evaluation may introduce too much overhead and degrade the overall performance.

Overall, these experiments show that effective use of \package{} requires matching the formulation, LMO implementation, Frank-Wolfe variant, tolerance settings, and heuristic strategy to 
the structure of the underlying problem.
}

\section{Discussion and conclusion}\label{sec:Discussion}

The Boscia framework brings a new modeling and solving paradigm to the MINLP solver landscape. 
In particular, when the linear constraints encode combinatorial structures which can be exploited,
it excels for large scale problems compared to outer approximation schemes. 
Furthermore, it is highly customizable giving the user full control over the optimization process.

In the following, we list some highlighted points and future work.
The objective function and its gradient are invoked frequently throughout the Frank-Wolfe iterations, making their efficient implementation crucial for overall performance.
The Frank-Wolfe iteration limit represents a fundamental trade-off: lower limits enable faster exploration of the branch-and-bound tree but may compromise the quality of lower bound improvements.
How much computational effort to spend at each node to tighten bounds is a challenging open question that will require both further empirical studies and theoretical understanding.

The FW variant selection presents another challenge, as a particular FW method that performs well on the root problem may not maintain its effectiveness on child nodes due to local geometric variations.
The framework aims to address this limitation in future work by detecting such situations and dynamically adjusting the variant selection.

Crucial polyhedral LMOs from \texttt{FrankWolfe.jl} are supported, and we aim to extend the support for as many other LMOs as possible in the future.

Furthermore, we are in the process of developing different modes for Boscia.
Currently available are the {\sc default} mode which is the one described here and a {\sc heuristic} mode.
The default values for the optional settings slightly differ between the two modes.

The current setup of the package assumes that the objective function and its gradient are deterministic and exact. 
This is a common assumption in the literature, but it is not always the case in applications.
In future work, we aim to investigate support for stochastic objective functions and their gradients.

The convexity requirement for the objective function represents another current limitation, with future plans to incorporate spatial branching for non-convex objectives.
Finally, relaxing requirements such as $L$-smoothness and differentiability would require complementary advances in the Frank-Wolfe methodology itself also planned as future work.

\section*{Acknowledgments}

Research reported in this paper was partially supported through the Research Campus Modal funded by the German Federal Ministry of Education 
and Research (fund numbers 05M14ZAM,05M20ZBM) and the Deutsche Forschungsgemeinschaft (DFG) through the DFG Cluster of Excellence MATH+ Project AA3-15.

%\clearpage
\bibliographystyle{icml2021}
\bibliography{references}

@article{hendrych2023convex,
      title={Convex mixed-integer optimization with {Frank-Wolfe} methods}, 
      author={Deborah Hendrych and Hannah Troppens and Mathieu Besançon and Sebastian Pokutta},
      journal={Mathematical Programming Computation},
      year={2025},
      primaryClass={math.OC}
}

@article{besanccon2025improved,
  title={Improved algorithms and novel applications of the {FrankWolfe.jl} library},
  author={Besan{\c{c}}on, Mathieu and Designolle, S{\'e}bastien and Halbey, Jannis and Hendrych, Deborah and Kuzinowicz, Dominik and Pokutta, Sebastian and Troppens, Hannah and Herrmannsdoerfer, Daniel Viladrich and Wirth, Elias},
  journal={ACM Transactions on Mathematical Software},
  year={2025}
}

@InProceedings{hendrych2025secant,
  title = 	 {Secant Line Search for {F}rank-{W}olfe Algorithms},
  author =       {Hendrych, Deborah and Pokutta, Sebastian and Besan\c{c}on, Mathieu and Mart\'{\i}nez-Rubio, David},
  booktitle = 	 {Proceedings of the 42nd International Conference on Machine Learning},
  pages = 	 {23005--23029},
  year = 	 {2025},
  volume = 	 {267},
  series = 	 {Proceedings of Machine Learning Research},
  month = 	 {13--19 Jul},
  publisher =    {PMLR},
}

@article{2025_MexiEtAl_Frankwolfeheuristic_2508-01299,
  archiveprefix = {arXiv},
  eprint = {2508.01299},
  primaryclass = {math.OC},
  year = {2025},
  author = {Mexi, Gioni and Hendrych, Deborah and Designolle, Sébastien and Besançon, Mathieu and Pokutta, Sebastian},
  title = {A {Frank-Wolfe}-based {P}rimal {H}euristic for {Q}uadratic {Mixed-Integer} {O}ptimization},
  date = {2025-08-02}
}

@inproceedings{2024_SharmaHendrychBesanconPokutta_NetworkdesignMicoFrankwolfe,
  year = {2024},
  booktitle = {Proceedings of the INFORMS Optimization Society Conference},
  archiveprefix = {arXiv},
  eprint = {2402.00166},
  primaryclass = {math.OC},
  author = {Sharma, Kartikey and Hendrych, Deborah and Besançon, Mathieu and Pokutta, Sebastian},
  title = {Network {D}esign for the {T}raffic {A}ssignment {P}roblem with {M}ixed-{I}nteger {Frank-Wolfe}},
  date = {2024-01-31}
}

@misc{transportnetlibrary,
    author = {{Transportation Networks for Research Core Team}},
    title = {Transportation Networks for Research},
    howpublished = {\url{https://github.com/bstabler/TransportationNetworks}},
    note = {Online; accessed 11 March 2023} ,
    year=2023,
}

@inproceedings{2023_HendrychBesanconPokutta_Optimalexperimentdesign,
  year = {2024},
  booktitle = {Proceedings of the Symposium on Experimental Algorithms},
  doi = {10.4230/LIPIcs.SEA.2024.16},
  archiveprefix = {arXiv},
  eprint = {2312.11200},
  primaryclass = {math.OC},
  author = {Hendrych, Deborah and Besançon, Mathieu and Pokutta, Sebastian},
  title = {Solving the {O}ptimal {E}xperiment {D}esign {P}roblem with {Mixed-Integer} {C}onvex {M}ethods},
  code = {https://github.com/ZIB-IOL/OptimalDesignWithBoscia}
}

@misc{BosciaDocumentation,
  title = {{Boscia.jl} {D}ocumentation},
  howpublished = {Available at: \url{https://zib-iol.github.io/Boscia.jl/stable/}},
  year = {2022},
  author={{IOL Lab}},
  url = "https://zib-iol.github.io/Boscia.jl/stable/"
}

@misc{BosciaGitHub,
  title = {{Boscia.jl} {G}itHub {R}epository},
  howpublished = {Available at: \url{https://github.com/ZIB-IOL/Boscia.jl}},
  year = {2022},
  author={{IOL Lab}},
  url = "https://github.com/ZIB-IOL/Boscia.jl"
}

@misc{FrankWolfeGitHub,
  title = {{FrankWolfe.jl} {G}itHub {R}epository},
  howpublished = {Available at: \url{https://github.com/ZIB-IOL/FrankWolfe.jl}},
  year = {2021},
  author={{IOL Lab}},
  url = "https://github.com/ZIB-IOL/FrankWolfe.jl"
}

@misc{CLOraclesGitHub,
  title = {{CombinatorialLinearOracles.jl} {G}itHub Repository},
  howpublished = {Available at: \url{https://github.com/ZIB-IOL/CombinatorialLinearOracles.jl}},
  year = {2025},
  author={{IOL Lab}},
  url = "https://github.com/ZIB-IOL/CombinatorialLinearOracles.jl"
}

@article{kronqvist2019review,
  title={A review and comparison of solvers for convex {MINLP}},
  author={Kronqvist, Jan and Bernal, David E and Lundell, Andreas and Grossmann, Ignacio E},
  journal={Optimization and Engineering},
  volume={20},
  pages={397--455},
  year={2019},
  publisher={Springer}
}

@article{bonami2008algorithmic,
  title={An algorithmic framework for convex mixed integer nonlinear programs},
  author={Bonami, Pierre and Biegler, Lorenz T and Conn, Andrew R and Cornu{\'e}jols, G{\'e}rard and Grossmann, Ignacio E and Laird, Carl D and Lee, Jon and Lodi, Andrea and Margot, Fran{\c{c}}ois and Sawaya, Nicolas and others},
  journal={Discrete optimization},
  volume={5},
  number={2},
  pages={186--204},
  year={2008},
  publisher={Elsevier}
}

@article{aps2022mosek,
  title={{MOSEK} {P}ortfolio {O}ptimization {C}ookbook},
  author={ApS, MOSEK},
  year={2022},
  url={https://docs.mosek.com/MOSEKModelingCookbook-letter.pdf}
}

@article{lundell2022supporting,
  title={The supporting hyperplane optimization toolkit for convex {MINLP}},
  author={Lundell, Andreas and Kronqvist, Jan and Westerlund, Tapio},
  journal={Journal of Global Optimization},
  volume={84},
  number={1},
  pages={1--41},
  year={2022},
  publisher={Springer}
}

@article{lundell2022polyhedral,
  title={Polyhedral approximation strategies for nonconvex mixed-integer nonlinear programming in {SHOT}},
  author={Lundell, Andreas and Kronqvist, Jan},
  journal={Journal of Global Optimization},
  volume={82},
  number={4},
  pages={863--896},
  year={2022},
  publisher={Springer}
}

@misc{gurobi,
  author = {{Gurobi Optimization, LLC}},
  title = {{Gurobi Optimizer Reference Manual}},
  year = 2024,
  url = "https://www.gurobi.com"
}

@article{wachter2006implementation,
  author       = {Andreas W{\"{a}}chter and
                  Lorenz T. Biegler},
  title        = {On the implementation of an interior-point filter line-search algorithm
                  for large-scale nonlinear programming},
  journal      = {Mathematical Programming},
  volume       = {106},
  number       = {1},
  pages        = {25--57},
  year         = {2006},
  publisher    = {Springer},
  url          = {https://doi.org/10.1007/s10107-004-0559-y},
  doi          = {10.1007/S10107-004-0559-Y},
  bibsource    = {dblp computer science bibliography, https://dblp.org}
}

@misc{wei2023outer,
      title={An Outer Approximation Method for Solving Mixed-Integer Convex Quadratic Programs with Indicators}, 
      author={Linchuan Wei and Simge Küçükyavuz},
      year={2023},
      eprint={2312.04812},
      archivePrefix={arXiv},
      primaryClass={math.OC},
      url={https://arxiv.org/abs/2312.04812}, 
}

@article{frank1956algorithm,
  title={An algorithm for quadratic programming},
  author={Frank, Marguerite and Wolfe, Philip},
  journal={Naval research logistics quarterly},
  volume={3},
  number={1-2},
  pages={95--110},
  year={1956},
  publisher={Wiley Online Library}
}

@article{levitin1966constrained,
  title={Constrained minimization methods},
  author={Levitin, Evgeny S and Polyak, Boris T},
  journal={USSR Computational mathematics and mathematical physics},
  volume={6},
  number={5},
  pages={1--50},
  year={1966},
  publisher={No longer published by Elsevier}
}

@article{wolfe1970convergence,
  title={Convergence theory in nonlinear programming},
  author={Wolfe, Philip},
  journal={Integer and nonlinear programming},
  pages={1--36},
  year={1970},
  publisher={North-Holland, Amsterdam}
}

@article{besanccon2022frankwolfe,
  title={{FrankWolfe}.jl: A {H}igh-{P}erformance and {F}lexible {T}oolbox for {F}rank--{W}olfe {A}lgorithms and {C}onditional {G}radients},
  author={Besan{\c{c}}on, Mathieu and Carderera, Alejandro and Pokutta, Sebastian},
  journal={INFORMS Journal on Computing},
  year={2022},
  publisher={INFORMS}
}

@inproceedings{braun2019blended,
  title={Blended conditonal gradients},
  author={Braun, G{\'a}bor and Pokutta, Sebastian and Tu, Dan and Wright, Stephen},
  booktitle={International Conference on Machine Learning},
  pages={735--743},
  year={2019},
  organization={PMLR}
}

@inproceedings{braun2017lazifying,
  title={Lazifying conditional gradient algorithms},
  author={Braun, G{\'a}bor and Pokutta, Sebastian and Zink, Daniel},
  booktitle={International conference on machine learning},
  pages={566--575},
  year={2017},
  organization={PMLR}
}

@article{garber2016linear,
  title={Linear-memory and decomposition-invariant linearly convergent conditional gradient algorithm for structured polytopes},
  author={Garber, Dan and Meshi, Ofer},
  journal={Advances in neural information processing systems},
  volume={29},
  year={2016}
}

@inproceedings{jaggi2013revisiting,
  title={Revisiting {F}rank--{W}olfe: {P}rojection-free sparse convex optimization},
  author={Jaggi, Martin},
  booktitle={International Conference on Machine Learning},
  pages={427--435},
  year={2013},
  organization={PMLR}
}

@book{braun2022conditional,
author = {Braun, G{\'a}bor and Carderera, Alejandro and Combettes, Cyrille W. and Hassani, Hamed and Karbasi, Amin and Mokhtari, Aryan and Pokutta, Sebastian},
title = {Conditional {G}radient {M}ethods},
publisher = {Society for Industrial and Applied Mathematics},
year = {2025},
doi = {10.1137/1.9781611978568},
address = {Philadelphia, PA},
edition   = {},
URL = {https://epubs.siam.org/doi/abs/10.1137/1.9781611978568},
eprint = {https://epubs.siam.org/doi/pdf/10.1137/1.9781611978568}
}

@article{sun2019generalized,
  title={Generalized self-concordant functions: a recipe for {N}ewton-type methods},
  author={Sun, Tianxiao and Tran-Dinh, Quoc},
  journal={Mathematical Programming},
  volume={178},
  number={1-2},
  pages={145--213},
  year={2019},
  publisher={Springer}
}

@inproceedings{
halbey2025efficient,
title={Efficient {Q}uadratic {C}orrections for {Frank-Wolfe} {A}lgorithms},
author={Jannis Halbey and Seta Rakotomandimby and Mathieu Besan{\c{c}}on and S{\'e}bastien Designolle and Sebastian Pokutta},
booktitle={The Thirty-ninth Annual Conference on Neural Information Processing Systems},
year={2025}
}

@article{holloway1974extension,
  title={An extension of the {Frank} and {Wolfe} method of feasible directions},
  author={Holloway, Charles A},
  journal={Mathematical Programming},
  volume={6},
  number={1},
  pages={14--27},
  year={1974},
  publisher={Springer}
}

@InProceedings{pedregosa2020linearly,
  title = 	 {Linearly {C}onvergent {F}rank-{W}olfe with {B}acktracking {L}ine-{S}earch},
  author =       {Pedregosa, Fabian and Negiar, Geoffrey and Askari, Armin and Jaggi, Martin},
  booktitle = 	 {Proceedings of the Twenty Third International Conference on Artificial Intelligence and Statistics},
  pages = 	 {1--10},
  year = 	 {2020},
  editor = 	 {Chiappa, Silvia and Calandra, Roberto},
  volume = 	 {108},
  series = 	 {Proceedings of Machine Learning Research},
  month = 	 {26--28 Aug},
  publisher =    {PMLR},
  }

@Article{pokutta2024frank,
  author  = {Pokutta, Sebastian},
  title   = {The {F}rank-{W}olfe {A}lgorithm: {A} {S}hort {I}ntroduction},
  journal = {Jahresbericht der Deutschen Mathematiker-Vereinigung},
  year    = {2024},
  volume  = {126},
  number  = {1},
  pages   = {3--35},
  doi     = {10.1365/s13291-023-00275-x},
}

@InProceedings{xiao2025GI,
author={Xiao, Wenjie and Besan{\c{c}}on, Mathieu and Gel{\ss}, Patrick and Hendrych, Deborah and Klus, Stefan and Pokutta, Sebastian},
title={Graph {I}somorphism: {M}ixed-{I}nteger {C}onvex {O}ptimization from {F}irst-{O}rder {M}ethods},
booktitle={Integration of Constraint Programming, Artificial Intelligence, and Operations Research},
year={2026},
publisher={Springer Nature Switzerland},
pages={614--630},
isbn={978-3-032-27242-3}
}

@article{hasan2017graphettes,
author = {Hasan, Adib and Chung, Po-Chien and Hayes, Wayne},
year = {2017},
month = {08},
pages = {},
title = {Graphettes: {C}onstant-time determination of graphlet and orbit identity including (possibly disconnected) graphlets up to size 8},
volume = {12},
journal = {PLOS ONE},
doi = {10.1371/journal.pone.0181570}
}

@article{klus2025continuous,
  title={Continuous optimization methods for the graph isomorphism problem},
  author={Klus, Stefan and Gel{\ss}, Patrick},
  journal={Information and Inference: A Journal of the IMA},
  volume={14},
  number={2},
  pages={iaaf011},
  year={2025},
  publisher={Oxford University Press}
}

@article{sagnol2015computing,
  title={Computing exact {D}-optimal designs by mixed integer second-order cone programming},
  author={Sagnol, Guillaume and Harman, Radoslav},
  journal={The Annals of Statistics},
  volume={43},
  number={5},
  pages={2198--2224},
  year={2015}
}

@article{sagnol2011computing,
  title={Computing optimal designs of multiresponse experiments reduces to second-order cone programming},
  author={Sagnol, Guillaume},
  journal={Journal of Statistical Planning and Inference},
  volume={141},
  number={5},
  pages={1684--1708},
  year={2011},
  publisher={Elsevier}
}

@article{ponte2023branch,
  title={Branch-and-bound for {D}-Optimality with fast local search and variable-bound tightening},
  author={Ponte, Gabriel and Fampa, Marcia and Lee, Jon},
  journal={Mathematical Programming},
  pages={337-374},
  year={2026},
  number={1},
  volume={218},
  publisher={Springer}
}

@inproceedings{pillai2025computing,
  title={Computing experiment-constrained d-optimal designs},
  author={Pillai, Aditya and Ponte, Gabriel and Fampa, Marcia and Lee, Jon and Singh, Mohit and Xie, Weijun},
  booktitle={2025 Proceedings of the Conference on Applied and Computational Discrete Algorithms (ACDA)},
  pages={182--195},
  year={2025},
  organization={SIAM}
}

@incollection{byrd2006knitro,
  title={KNITRO: An integrated package for nonlinear optimization},
  author={Byrd, Richard H and Nocedal, Jorge and Waltz, Richard A},
  booktitle={Large-scale nonlinear optimization},
  pages={35--59},
  year={2006},
  publisher={Springer}
}

@article{coey2022solving,
    title={Solving natural conic formulations with {H}ypatia.jl},
    author={Chris Coey and Lea Kapelevich and Juan Pablo Vielma},
    year={2022},
    journal={INFORMS Journal on Computing},
    publisher={INFORMS},
    volume={34},
    number={5},
    pages={2686--2699},
    doi={https://doi.org/10.1287/ijoc.2022.1202}
}

@article{CoeyLubinVielma2020,
    title={Outer approximation with conic certificates for mixed-integer convex problems},
    author={Coey, Chris and Lubin, Miles and Vielma, Juan Pablo},
    journal={Mathematical Programming Computation},
    volume={12},
    number={2},
    pages={249--293},
    year={2020},
    publisher={Springer}
}

@article{Lubin2023,
    author = {Miles Lubin and Oscar Dowson and Joaquim {Dias Garcia} and Joey Huchette and Beno{\^i}t Legat and Juan Pablo Vielma},
    title = {{JuMP} 1.0: {R}ecent improvements to a modeling language for mathematical optimization},
    journal = {Mathematical Programming Computation},
    year = {2023},
    doi = {10.1007/s12532-023-00239-3}
}

@techreport{BolusaniEtal2024OO,
  author = {Suresh Bolusani and Mathieu Besan{\c{c}}on and Ksenia Bestuzheva and Antonia Chmiela and Jo{\~{a}}o Dion{\'{i}}sio and Tim Donkiewicz and Jasper van Doornmalen and Leon Eifler and Mohammed Ghannam and Ambros Gleixner and Christoph Graczyk and Katrin Halbig and Ivo Hedtke and Alexander Hoen and Christopher Hojny and Rolf van der Hulst and Dominik Kamp and Thorsten Koch and Kevin Kofler and Jurgen Lentz and Julian Manns and Gioni Mexi and Erik~M\"{u}hmer and Marc E. Pfetsch and Franziska Schl{\"o}sser and Felipe Serrano and Yuji Shinano and Mark Turner and Stefan Vigerske and Dieter Weninger and Lixing Xu},
  title = {{The SCIP Optimization Suite 9.0}},
  type = {Technical Report},
  institution = {Optimization Online},
  month = {February},
  year = {2024},
  url = {https://optimization-online.org/2024/02/the-scip-optimization-suite-9-0/}
}

@article{huangfu2018parallelizing,
  title={Parallelizing the dual revised simplex method},
  author={Huangfu, Qi and Hall, Julian},
  journal={Mathematical Programming Computation},
  volume={10},
  number={1},
  pages={119--142},
  year={2018},
  publisher={Springer}
}

@misc{Xpress,
  author = {FICO},
  title = {{FICO} {Xpress} Optimizer},
  howpublished = {Available at: \url{http://www.fico.com/en/Products/DMTools/xpress-overview/Pages/Xpress-Optimizer.aspx}}
}

@article{bonami2013branching,
  author       = {Pierre Bonami and
                  Jon Lee and
                  Sven Leyffer and
                  Andreas W{\"{a}}chter},
  title        = {On branching rules for convex mixed-integer nonlinear optimization},
  journal      = {{ACM} Journal of Experimental Algorithmics},
  volume       = {18},
  year         = {2013},
  url          = {https://doi.org/10.1145/2532568},
  doi          = {10.1145/2532568},
  bibsource    = {dblp computer science bibliography, https://dblp.org}
}

@inproceedings{applegate1998solution,
  title={On the solution of traveling salesman problems},
  author={Applegate, David and Bixby, Robert and Chv{\'a}tal, Va{\v{s}}ek and Cook, William},
  booktitle={Proceedings of the International Congress of Mathematicians 1998},
  pages={645--656},
  year={1998},
  organization={European Mathematical Society-EMS-Publishing House GmbH}
}

@article{benichou1971experiments,
  author       = {Michel B{\'{e}}nichou and
                  Jean{-}Michel Gauthier and
                  Paul Girodet and
                  Gerard Hentges and
                  Gerard Ribi{\`{e}}re and
                  O. Vincent},
  title        = {Experiments in mixed-integer linear programming},
  journal      = {Math. Program.},
  volume       = {1},
  number       = {1},
  pages        = {76--94},
  year         = {1971},
  url          = {https://doi.org/10.1007/BF01584074},
  doi          = {10.1007/BF01584074},
  bibsource    = {dblp computer science bibliography, https://dblp.org}
}

@article{belotti2009branching,
  title={Branching and bounds tightening techniques for non-convex {MINLP}},
  author={Belotti, Pietro and Lee, Jon and Liberti, Leo and Margot, Fran{\c{c}}ois and W{\"a}chter, Andreas},
  journal={Optimization Methods \& Software},
  volume={24},
  number={4-5},
  pages={597--634},
  year={2009},
  publisher={Taylor \& Francis}
}

@inproceedings{achterberg2009hybrid,
  author       = {Tobias Achterberg and
                  Timo Berthold},
  editor       = {Willem Jan van Hoeve and
                  John N. Hooker},
  title        = {Hybrid {B}ranching},
  booktitle    = {Integration of {AI} and {OR} {T}echniques in {C}onstraint {P}rogramming
                  for {C}ombinatorial {O}ptimization {P}roblems, 6th International Conference,
                  {CPAIOR} 2009, Pittsburgh, PA, USA, May 27-31, 2009, Proceedings},
  series       = {Lecture Notes in Computer Science},
  pages        = {309--311},
  publisher    = {Springer},
  year         = {2009},
  url          = {https://doi.org/10.1007/978-3-642-01929-6\_23},
  doi          = {10.1007/978-3-642-01929-6\_23},
  bibsource    = {dblp computer science bibliography, https://dblp.org}
}

@article{gally2018framework,
  author       = {Tristan Gally and
                  Marc E. Pfetsch and
                  Stefan Ulbrich},
  title        = {A framework for solving mixed-integer semidefinite programs},
  journal      = {Optimization Methods and Software},
  volume       = {33},
  number       = {3},
  pages        = {594--632},
  year         = {2018},
  url          = {https://doi.org/10.1080/10556788.2017.1322081},
  doi          = {10.1080/10556788.2017.1322081},
  bibsource    = {dblp computer science bibliography, https://dblp.org}
}

@misc{SCIPSDPGithub,
  title = {{SCIPSDP} {G}itHub Repository},
  howpublished = {Available at: \url{https://github.com/scipopt/SCIP-SDP}},
  year = {2022},
  url = "https://github.com/scipopt/SCIP-SDP"
}

@InProceedings{lamperski2025simple,
  title = 	 {Simple Randomized Rounding for {M}ax-{M}in Eigenvalue Augmentation},
  author =       {Lamperski, Jourdain and Yang, Haeseong and Prokopyev, Oleg},
  booktitle = 	 {Proceedings of the 42nd International Conference on Machine Learning},
  pages = 	 {32406--32419},
  year = 	 {2025},
  editor = 	 {Singh, Aarti and Fazel, Maryam and Hsu, Daniel and Lacoste-Julien, Simon and Berkenkamp, Felix and Maharaj, Tegan and Wagstaff, Kiri and Zhu, Jerry},
  volume = 	 {267},
  series = 	 {Proceedings of Machine Learning Research},
  month = 	 {13--19 Jul},
  publisher =    {PMLR},
  url = 	 {https://proceedings.mlr.press/v267/lamperski25a.html}
}

@book{fedorov1972theory,
  title={Theory of {O}ptimal {E}xperiments},
  author={Fedorov, Valerii Vadimovich},
  year={1972},
  publisher={Academic Press},
  address={New York},
  translator={Studden, W. J. and Klimko, E. M.}
}

@article{mathoptinterface,
    title={{MathOptInterface}: a data structure for mathematical optimization problems},
    author={Legat, Beno{\^\i}t and Dowson, Oscar and Dias Garcia, Joaquim and Lubin, Miles},
    journal={INFORMS Journal on Computing},
    year={2021},
    volume={34},
    number={2},
    pages={672--689},
    doi={10.1287/ijoc.2021.1067},
    publisher={INFORMS}
}

@article{kiefer1953sequential,
 ISSN = {00029939, 10886826},
 URL = {http://www.jstor.org/stable/2032161},
 author = {J. Kiefer},
 journal = {Proceedings of the American Mathematical Society},
 number = {3},
 pages = {502--506},
 publisher = {American Mathematical Society},
 title = {Sequential Minimax Search for a Maximum},
 urldate = {2026-09-09},
 volume = {4},
 year = {1953}
}

@article{Avriel01101966,
author = {Mordecai Avriel and Douglass J. Wilde},
title = {Optimality Proof for the Symmetric Fibonacci Search Technique},
journal = {The Fibonacci Quarterly},
volume = {4},
number = {3},
pages = {265--269},
year = {1966},
publisher = {Taylor \& Francis},
doi = {10.1080/00150517.1966.12431364},
URL = {https://doi.org/10.1080/00150517.1966.12431364},
}

\clearpage
\appendix
\section{\revision{Network design problem}}
\label{app:network_design}

\lstinline{bounded_compute_extreme_point} receives the linear cost vector \lstinline{d}, the lower and upper bounds of the current node, and the indices of the integer variables, and returns the linear minimizer \lstinline{v}.

\noindent
\lstinline{is_simple_linear_feasible} checks whether a given point satisfies the continuous relaxations of the constraints.

\noindent
\begin{minipage}{\linewidth}
\begingroup
\renewcommand{\jlbasicfont}{\ttfamily\small\selectfont}
\begin{jllisting}[language=julia, style=jlcodestyle]{}
struct NetworkDesignLMO <: FrankWolfe.LinearMinimizationOracle 
    ...
end
function bounded_compute_extreme_point(lmo::NetworkDesignLMO, d, lb,
                     ub, int_vars; kwargs...)
    ...
    return v
end

function is_simple_linear_feasible(lmo::NetworkDesignLMO, v)
    ... 
    return true
end
\end{jllisting}
\endgroup
\vspace{-20pt}
\captionsetup{skip=2pt}
\captionof{lstlisting}{Implementation interfaces for the network design LMO}
\label{lst:network_design_lmo}
\end{minipage}

\vspace{10pt}

\section{\revision{Graph isomorphism problem}}
\label{app:graph_isomorphism}

\subsection{Extreme-Point Computation}

\lstinline{compute_extreme_point} reduces the problem to the active rows and columns, solves the resulting assignment problem with the Hungarian algorithm, 
and then lifts the solution back to the full matrix.

\begin{algorithm}[htbp]
    \caption{Implementation of \lstinline{compute_extreme_point} via the Hungarian algorithm}
    \label{alg:extreme_point}
    \begin{algorithmic}[1]
        \Require \texttt{lmo} encoding the feasible region and direction matrix $D$
        \State $n_{\mathrm{reduced}} \gets |\texttt{lmo.index\_map\_rows}|$
        \State $D_{\mathrm{reduced}} \gets \mathrm{zeros}(n_{\mathrm{reduced}}, n_{\mathrm{reduced}})$
        \ForAll{$i,j \in \{1,\dots,n_{\mathrm{reduced}}\}$} \Comment{Build reduced direction matrix}
            \State $\mathrm{linear\_index} \gets \mathrm{map\_to\_original\_linear\_index}(i,j)$
            \State $i_{\mathrm{orig}} \gets \texttt{lmo.index\_map\_rows}[i]$
            \State $j_{\mathrm{orig}} \gets \texttt{lmo.index\_map\_cols}[j]$
            \If{$\texttt{lmo.upper\_bounds}[\mathrm{linear\_index}] = 0$}
                \State $D_{\mathrm{reduced}}[i,j] \gets \infty$
            \Else
                \State $D_{\mathrm{reduced}}[i,j] \gets D[i_{\mathrm{orig}}, j_{\mathrm{orig}}]$
            \EndIf
        \EndFor
        \State $M_{\mathrm{reduced}} \gets \textsc{Hungarian}(D_{\mathrm{reduced}})$ \Comment{Solve assignment problem using Hungarian algorithm}
        \State $M \gets \mathrm{zeros}(n, n)$ \Comment{Recover the reduced solution}
        \For{$k = 1,\dots,|\texttt{lmo.fixed\_to\_one\_rows}|$}
            \State $i, j \gets \texttt{lmo.fixed\_to\_one\_rows}[k], \texttt{lmo.fixed\_to\_one\_cols}[k]$
            \State $M[i, j] \gets 1$
        \EndFor
        \ForAll{$i,j \in \{1,\dots,n_{\mathrm{reduced}}\}$}
            \State $i_{\mathrm{orig}}, j_{\mathrm{orig}} \gets \texttt{lmo.index\_map\_rows}[i], \texttt{lmo.index\_map\_cols}[j]$
            \State $M[i_{\mathrm{orig}}, j_{\mathrm{orig}}] \gets M_{\mathrm{reduced}}[i, j]$
        \EndFor
        \State \Return $M$
    \end{algorithmic}
\end{algorithm}

\subsection{Maximum-Step Computation}
\label{app:maximum_step_computation}

\lstinline{dicg_maximum_step} checks how far one can move along the direction $D$ from the current point $X$ before hitting a bound, and 
returns zero if no positive step is feasible.

\begin{minipage}{\linewidth}
\begingroup
\renewcommand{\jlbasicfont}{\ttfamily\small\selectfont}
\begin{jllisting}[language=julia, style=jlcodestyle]{}
function dicg_maximum_step(lmo::BirkhoffLMO, D, X; kwargs...)
    n = lmo.dim
    T = promote_type(eltype(X), eltype(D))
    gamma_max = one(T)
    for idx in eachindex(X)
        if D[idx] != 0
            # iterate already on the boundary
            if (D[idx] < 0 && X[idx] ≈ 1) ||
               (D[idx] > 0 && X[idx] ≈ 0)
                return zero(gamma_max)
            end
            # clipping with the zero boundary
            if D[idx] > 0
                gamma_max = min(gamma_max, X[idx] / D[idx])
            else
                gamma_max = min(gamma_max, -(1 - X[idx]) / D[idx])
            end
        end
    end
    return gamma_max
end
\end{jllisting}
\endgroup
\vspace{-20pt}
\captionsetup{skip=2pt}
\captionof{lstlisting}{Implementation of \lstinline{dicg_maximum_step} for the Birkhoff polytope}
\label{lst:dicg_maximum_step}
\end{minipage}

\vspace{10pt}
\subsection{Implementation of Bound Management}

\noindent
\lstinline{build_global_bounds} builds the global bounds from the LMO and integer variables at initialization, stores them at the tree level, 
and reuses them when constructing node LMOs.

\begin{minipage}{\linewidth}
\begingroup
\renewcommand{\jlbasicfont}{\ttfamily\small\selectfont}
\begin{jllisting}[language=julia, style=jlcodestyle]{}
# Read global bounds from the problem.
function build_global_bounds(lmo::BirkhoffLMO, integer_variables)
    global_bounds = Boscia.IntegerBounds()
    for (idx, int_var) in enumerate(lmo.int_vars)
        push!(global_bounds, (int_var, lmo.lower_bounds[idx]), :greaterthan)
        push!(global_bounds, (int_var, lmo.upper_bounds[idx]), :lessthan)
    end
    return global_bounds
end
\end{jllisting}
\endgroup
\vspace{-20pt}
\captionsetup{skip=2pt}
\captionof{lstlisting}{Implementation of \lstinline{build_global_bounds}}
\label{lst:build_global_bounds}
\end{minipage}

\vspace{10pt}

\lstinline{set_bound!} receives the LMO, the variable index, the bound value, and the constraint sense, and sets the corresponding bound. 
It also updates \lstinline{fixed_to_one_rows} and \lstinline{fixed_to_one_cols}.

\begin{minipage}{\linewidth}
\begingroup
\renewcommand{\jlbasicfont}{\ttfamily\small\selectfont}
\begin{jllisting}[language=julia, style=jlcodestyle]{}
function set_bound!(lmo::BirkhoffLMO, c_idx, value, sense::Symbol)
    # Reset the lmo if necessary
    if lmo.updated_lmo
        empty!(lmo.fixed_to_one_rows)
        empty!(lmo.fixed_to_one_cols)
        lmo.updated_lmo = false
    end
    if sense == :greaterthan
        lmo.lower_bounds[c_idx] = value
        if value == 1
            n0 = lmo.dim
            fixed_int_var = lmo.int_vars[c_idx]
            # Convert linear index to (row, col) based on storage format
            j = ceil(Int, fixed_int_var / n0)  # column index
            i = Int(fixed_int_var - n0 * (j - 1))  # row index
            push!(lmo.fixed_to_one_rows, i)
            push!(lmo.fixed_to_one_cols, j)
        end
    elseif sense == :lessthan
        lmo.upper_bounds[c_idx] = value
    else
        error("Allowed values for sense are :lessthan and :greaterthan.")
    end
end
\end{jllisting}
\endgroup
\vspace{-20pt}
\captionsetup{skip=2pt}
\captionof{lstlisting}{Implementation of \lstinline{set_bound!}}
\label{lst:set_bound!}
\end{minipage}

\vspace{10pt}

\lstinline{delete_bounds!} removes selected bounds from the LMO, rebuilds the row and column index maps after accounting for entries fixed to one, 
and marks the LMO as updated.

\begin{minipage}{\linewidth}
\begingroup
\renewcommand{\jlbasicfont}{\ttfamily\small\selectfont}
\begin{jllisting}[language=julia, style=jlcodestyle]{}
function delete_bounds!(lmo::BirkhoffLMO, cons_delete)
    for (d_idx, sense) in cons_delete
        if sense == :greaterthan
            lmo.lower_bounds[d_idx] = 0
        else
            lmo.upper_bounds[d_idx] = 1
        end
    end

    nfixed = length(lmo.fixed_to_one_rows)
    nreduced = lmo.dim - nfixed

    # Store the indices of the original matrix that 
    # are still in the reduced matrix
    index_map_rows = fill(1, nreduced)
    index_map_cols = fill(1, nreduced)
    idx_in_map_row = 1
    idx_in_map_col = 1
    for orig_idx in 1:lmo.dim
        if orig_idx ∉ lmo.fixed_to_one_rows
            index_map_rows[idx_in_map_row] = orig_idx
            idx_in_map_row += 1
        end
        if orig_idx ∉ lmo.fixed_to_one_cols
            index_map_cols[idx_in_map_col] = orig_idx
            idx_in_map_col += 1
        end
    end

    empty!(lmo.index_map_rows)
    empty!(lmo.index_map_cols)
    append!(lmo.index_map_rows, index_map_rows)
    append!(lmo.index_map_cols, index_map_cols)
    lmo.updated_lmo = true
    return true
end
\end{jllisting}
\endgroup
\vspace{-20pt}
\captionsetup{skip=2pt}
\captionof{lstlisting}{Implementation of \lstinline{delete_bounds!}}
\label{lst:delete_bounds!}
\end{minipage}

\vspace{10pt}

\subsection{Implementation of Tree Callback}

\lstinline{build_tree_callback} verifies whether the current incumbent objective value is approximately zero. 
If this condition is satisfied, the solving stage is set to \lstinline{USER_STOP} , and the solution process terminates.

\noindent
\begin{minipage}{\linewidth}
\begingroup
\renewcommand{\jlbasicfont}{\ttfamily\small\selectfont}
\begin{jllisting}[language=julia, style=jlcodestyle]{}
function build_tree_callback()
    return function (
        tree, 
        node; 
        worse_than_incumbent=false, 
        node_infeasible=false, 
        lb_update=false,
        )
        if tree.incumbent ≈ 0
            tree.root.problem.solving_stage = Boscia.USER_STOP
            println("Optimal solution found.")
        end
        if Boscia.tree_lb(tree) > eps()
            tree.root.problem.solving_stage = Boscia.USER_STOP
            println("Tree lower bound already positive. No solution possible.")
        end
    end
end
\end{jllisting}
\endgroup
\vspace{-20pt}
\captionsetup{skip=2pt}
\captionof{lstlisting}{Implementation of \lstinline{build_tree_callback}}
\label{lst:build_tree_callback}
\end{minipage}

\vspace{10pt}

\lstinline{build_branch_callback} checks whether the node lower bound is already positive. 
If this is the case, it prints a message and returns \lstinline{false} to prevent further branching.

\noindent
\begin{minipage}{\linewidth}
\begingroup
\renewcommand{\jlbasicfont}{\ttfamily\small\selectfont}
\begin{jllisting}[language=julia, style=jlcodestyle]{}
function build_branch_callback()
    return function (tree, node, vidx)
        x = Bonobo.get_relaxed_values(tree, node)
        primal = tree.root.problem.f(x)
        lower_bound = primal - node.dual_gap
        if lower_bound > eps()
            println("No need to branch here. 
                        Node lower bound already positive.")
        end
        return lower_bound <= eps()
    end
end
\end{jllisting}
\endgroup
\vspace{-20pt}
\captionsetup{skip=2pt}
\captionof{lstlisting}{Implementation of \lstinline{build_branch_callback}}
\label{lst:build_branch_callback}
\end{minipage}

\clearpage
\section{\revision{Optimal experiment design problem}}
\label{app:optimal_design_problem}

\subsection{Domain oracle}

\lstinline{domain_oracle} checks whether the information matrix $X(\vx) = A^\top \mathrm{diag}(\vx)\, A$ is positive definite, i.e., whether $\vx$ lies in the domain of the objective function.

\noindent
\begin{minipage}{\linewidth}
\begingroup
\renewcommand{\jlbasicfont}{\ttfamily\small\selectfont}
\begin{jllisting}[language=julia, style=jlcodestyle]{}
function domain_oracle(x)
    X = A' * diagm(x) * A
    X = Symmetric(X)
    return LinearAlgebra.isposdef(X)
end
\end{jllisting}
\endgroup
\vspace{-20pt}
\captionsetup{skip=2pt}
\captionof{lstlisting}{Implementation of \lstinline{domain_oracle}}
\label{lst:domain_oracle}
\end{minipage}

\vspace{10pt}
\subsection{Domain point generation}

\lstinline{domain_point} first selects $n$ linearly independent rows of $A$ that respect the upper bounds, then iteratively increases entries until the budget $\sum_i x_i = N$ is met.
If no domain-feasible point exists, it returns \lstinline{nothing}.

\noindent
\begin{minipage}{\linewidth}
\begingroup
\renewcommand{\jlbasicfont}{\ttfamily\small\selectfont}
\begin{jllisting}[language=julia, style=jlcodestyle]{}
function domain_point(local_bounds)
    ...
    return x
end
\end{jllisting}
\endgroup
\vspace{-20pt}
\captionsetup{skip=2pt}
\captionof{lstlisting}{Implementation interface of \lstinline{domain_point}}
\label{lst:domain_point}
\end{minipage}

\vspace{10pt}
\begin{algorithm}[htbp]
    \caption{Domain point generation for OEDP}
    \label{alg:domain_point}
    \begin{algorithmic}[1]
    \Require Local bounds from branching, matrix $A \in \mathbb{R}^{m \times n}$, global upper bounds $\vu \in \mathbb{R}^m$, budget $N$
    \Ensure Domain-feasible point $\vx \in \mathbb{R}^m$ or \texttt{nothing} if infeasible
    
    \State Initialize $\vx \gets \mathbf{0}$, $\text{lb} \gets \mathbf{0}$, $\text{ub} \gets \vu$
    \State Apply local bounds from branching to lb and ub
    \If{$\sum \text{lb}_i > N$ or $\neg \text{domain\_oracle}(\text{ub})$}
        \State \Return \texttt{nothing} \Comment{Node infeasible}
    \EndIf
    
    \State $\vx \gets \text{lb}$
    \State $S \gets$ \Call{LinearlyIndependentRows}{$A$, $\left(\text{ub}_i > 0\right)_{i=1}^m$} \Comment{Select $n$ linearly independent rows of $A$}
    
    \While{$\sum x_i \leq N$}
        \If{$\sum x_i = N$}
            \State \Return $\vx$ if $\text{domain\_oracle}(\vx)$, otherwise \texttt{nothing}
        \EndIf
        \If{$\vx[S] \neq \text{ub}[S]$}
            \State $\vx[S] \gets$ \Call{AddToMin}{$\vx[S]$, $\text{ub}[S]$} \Comment{Increase selected experiments}
        \Else
            \State $\vx \gets$ \Call{AddToMin}{$\vx$, $\text{ub}$} \Comment{Add a new experiment}
        \EndIf
    \EndWhile
    \State \Return $\vx$
    \end{algorithmic}
\end{algorithm}

\subsection{Starting-point projection}

Given a domain-feasible point~$\vx_0$, we compute a starting point by solving a Frank–Wolfe projection problem, 
using a custom callback to move the iterate a few steps into the domain.

\noindent
\begin{minipage}{\linewidth}
\begingroup
\renewcommand{\jlbasicfont}{\ttfamily\small\selectfont}
\begin{jllisting}[language=julia, style=jlcodestyle]{}
# Build initial start point using domain_point function
initial_bounds = Boscia.IntegerBounds(zeros(m), u, collect(1:m))
x0 = domain_point(initial_bounds)

# Solve auxiliary problem to find feasible active set
f_help(x) = 1/2 * LinearAlgebra.norm(x - x0)^2
grad_help!(storage, x) = storage .= x - x0
v0 = compute_extreme_point(lmo, collect(1.0:m))

# Custom callback to ensure domain feasibility
function build_inner_callback()
    domain_counter = 0
    return function inner_callback(state, active_set, kwargs...)
        if domain_oracle(state.x)
            if domain_counter > 5
                return false
            end
            domain_counter += 1
        end
    end
end

inner_callback = build_inner_callback()

x, _, _, _, _, active_set = FrankWolfe.blended_pairwise_conditional_gradient(
    f_help,
    grad_help!,
    lmo,
    v0,
    callback=inner_callback,
    lazy=true,
)
\end{jllisting}
\endgroup
\vspace{-20pt}
\captionsetup{skip=2pt}
\captionof{lstlisting}{Projection-based starting-point generation}
\label{lst:oed_starting_point}
\end{minipage}

\subsection{Results for different Frank-Wolfe tolerances}
\label{app:oed_fw_tolerances}
\begin{table}[htbp]
    \centering
    \resizebox{\linewidth}{!}{%
    \begin{tabular}{ccrrrrrrrrrrrr}
    \toprule
    \# exp. & \# par. & \multicolumn{2}{c}{5e-2} & \multicolumn{2}{c}{1e-2} & \multicolumn{2}{c}{5e-3} & \multicolumn{2}{c}{1e-3} & \multicolumn{2}{c}{5e-4} & \multicolumn{2}{c}{1e-4} \\
    \cmidrule(lr){3-4}\cmidrule(lr){5-6}\cmidrule(lr){7-8}\cmidrule(lr){9-10}\cmidrule(lr){11-12}\cmidrule(lr){13-14}
     &  & \makecell{\%\\solved} & Abs.\ gap & \makecell{\%\\solved} & Abs.\ gap & \makecell{\%\\solved} & Abs.\ gap & \makecell{\%\\solved} & Abs.\ gap & \makecell{\%\\solved} & Abs.\ gap & \makecell{\%\\solved} & Abs.\ gap \\
     \midrule
     50 & 5 & \textbf{100\%} & 6.68e-07 & \textbf{100\%} & 7.22e-07 & \textbf{100\%} & 6.32e-07 & \textbf{100\%} & 2.54e-07 & \textbf{100\%} & \textbf{1.99e-07} & \textbf{100\%} & 8.39e-07 \\
     50 & 12 & 0\% & 4.21e-04 & \textbf{80\%} & 3.84e-05 & \textbf{80\%} & 2.26e-05 & \textbf{80\%} & 7.66e-06 & \textbf{80\%} & \textbf{6.32e-06} & \textbf{80\%} & 9.58e-06 \\
     60 & 6 & \textbf{100\%} & 9.96e-07 & \textbf{100\%} & 9.84e-07 & \textbf{100\%} & 9.69e-07 & \textbf{100\%} & 9.14e-07 & \textbf{100\%} & \textbf{7.72e-07} & \textbf{100\%} & 9.00e-07 \\
     60 & 15 & 0\% & 8.79e-04 & 0\% & 2.01e-04 & 0\% & 1.20e-04 & \textbf{40\%} & 6.38e-05 & \textbf{40\%} & \textbf{5.29e-05} & 20\% & 6.62e-05 \\
     80 & 8 & 0\% & 2.79e-04 & \textbf{100\%} & 9.97e-07 & \textbf{100\%} & 9.91e-07 & \textbf{100\%} & 9.90e-07 & \textbf{100\%} & 9.86e-07 & \textbf{100\%} & \textbf{9.76e-07} \\
     80 & 20 & 0\% & 1.28e-03 & 0\% & 3.93e-04 & 0\% & 2.66e-04 & 0\% & 1.51e-04 & 0\% & \textbf{1.38e-04} & 0\% & 1.40e-04 \\
     100 & 10 & 0\% & 5.37e-04 & 0\% & 1.37e-04 & 0\% & 6.11e-05 & 20\% & 7.91e-06 & 60\% & 3.03e-06 & \textbf{100\%} & \textbf{9.99e-07} \\
     100 & 25 & 0\% & 1.41e-03 & 0\% & 4.92e-04 & 0\% & 3.46e-04 & 0\% & 2.31e-04 & 0\% & 2.08e-04 & 0\% & \textbf{2.04e-04} \\
     120 & 12 & 0\% & 6.24e-04 & 0\% & 1.80e-04 & 0\% & 9.70e-05 & 0\% & 2.12e-05 & 0\% & 1.40e-05 & 0\% & \textbf{7.48e-06} \\
     120 & 30 & 0\% & 1.44e-03 & 0\% & 5.23e-04 & 0\% & 3.52e-04 & 0\% & 1.94e-04 & 0\% & 1.74e-04 & 0\% & \textbf{1.57e-04} \\
     \bottomrule
    \end{tabular}%
    }
    \caption{\protect\revII{Abs.\ gap is the mean absolute gap for different Frank-Wolfe tolerances on A-optimal design with correlated data.
    The best result in each row is shown in bold.}}
    \label{tab:final_abs_gap_a_correlated}
\end{table}

\begin{table}[htbp]
    \centering
    \resizebox{\linewidth}{!}{%
    \begin{tabular}{ccrrrrrrrrrrrr}
    \toprule
    \# exp. & \# par. & \multicolumn{2}{c}{5e-2} & \multicolumn{2}{c}{1e-2} & \multicolumn{2}{c}{5e-3} & \multicolumn{2}{c}{1e-3} & \multicolumn{2}{c}{5e-4} & \multicolumn{2}{c}{1e-4} \\
    \cmidrule(lr){3-4}\cmidrule(lr){5-6}\cmidrule(lr){7-8}\cmidrule(lr){9-10}\cmidrule(lr){11-12}\cmidrule(lr){13-14}
     &  & \makecell{\%\\solved} & Abs.\ gap & \makecell{\%\\solved} & Abs.\ gap & \makecell{\%\\solved} & Abs.\ gap & \makecell{\%\\solved} & Abs.\ gap & \makecell{\%\\solved} & Abs.\ gap & \makecell{\%\\solved} & Abs.\ gap \\
     \midrule
     50 & 5 & \textbf{100\%} & \textbf{0.00e+00} & \textbf{100\%} & 2.65e-05 & \textbf{100\%} & 1.21e-05 & \textbf{100\%} & 1.32e-07 & \textbf{100\%} & 1.32e-07 & \textbf{100\%} & 1.32e-07 \\
     50 & 12 & \textbf{60\%} & \textbf{2.75e-04} & \textbf{60\%} & 2.95e-04 & \textbf{60\%} & 3.80e-04 & \textbf{60\%} & 7.82e-04 & \textbf{60\%} & 9.39e-04 & \textbf{60\%} & 1.20e-03 \\
     60 & 6 & \textbf{100\%} & \textbf{2.41e-08} & \textbf{100\%} & 2.41e-06 & \textbf{100\%} & \textbf{2.41e-08} & \textbf{100\%} & 4.70e-08 & \textbf{100\%} & 1.80e-07 & \textbf{100\%} & 3.43e-08 \\
     60 & 15 & 0\% & 6.60e-03 & 0\% & \textbf{5.98e-03} & 0\% & 6.59e-03 & 0\% & 7.43e-03 & 0\% & 7.89e-03 & 0\% & 8.32e-03 \\
     80 & 8 & \textbf{100\%} & 4.16e-07 & \textbf{100\%} & 4.99e-07 & \textbf{100\%} & \textbf{5.33e-08} & \textbf{100\%} & 1.38e-07 & \textbf{100\%} & 1.97e-07 & \textbf{100\%} & 3.34e-07 \\
     80 & 20 & 0\% & \textbf{1.17e-02} & 0\% & 1.20e-02 & 0\% & \textbf{1.17e-02} & 0\% & 1.28e-02 & 0\% & 1.34e-02 & 0\% & 1.34e-02 \\
     100 & 10 & 40\% & 7.35e-04 & 40\% & 3.03e-04 & \textbf{60\%} & \textbf{2.82e-04} & 20\% & 4.73e-04 & 20\% & 5.58e-04 & 20\% & 8.26e-04 \\
     100 & 25 & 0\% & 1.92e-02 & 0\% & \textbf{1.74e-02} & 0\% & 1.92e-02 & 0\% & 2.02e-02 & 0\% & 2.05e-02 & 0\% & 2.33e-02 \\
     120 & 12 & 0\% & 2.95e-03 & 0\% & 2.12e-03 & 0\% & \textbf{2.00e-03} & 0\% & 2.08e-03 & 0\% & 2.17e-03 & 0\% & 2.43e-03 \\
     120 & 30 & 0\% & 2.41e-02 & 0\% & 2.46e-02 & 0\% & \textbf{2.39e-02} & 0\% & 2.86e-02 & 0\% & 2.71e-02 & 0\% & 3.09e-02 \\
     \bottomrule
    \end{tabular}%
    }
    \caption{\protect\revII{Abs.\ gap is the mean absolute gap for different Frank-Wolfe tolerances on A-optimal design with independent data.
    The best result in each row is shown in bold.}}
    \label{tab:final_abs_gap_a_independent}
\end{table}

\begin{table}[htbp]
    \centering
    \resizebox{\linewidth}{!}{%
    \begin{tabular}{ccrrrrrrrrrrrr}
    \toprule
    \# exp. & \# par. & \multicolumn{2}{c}{5e-2} & \multicolumn{2}{c}{1e-2} & \multicolumn{2}{c}{5e-3} & \multicolumn{2}{c}{1e-3} & \multicolumn{2}{c}{5e-4} & \multicolumn{2}{c}{1e-4} \\
    \cmidrule(lr){3-4}\cmidrule(lr){5-6}\cmidrule(lr){7-8}\cmidrule(lr){9-10}\cmidrule(lr){11-12}\cmidrule(lr){13-14}
     &  & \makecell{\%\\solved} & Abs.\ gap & \makecell{\%\\solved} & Abs.\ gap & \makecell{\%\\solved} & Abs.\ gap & \makecell{\%\\solved} & Abs.\ gap & \makecell{\%\\solved} & Abs.\ gap & \makecell{\%\\solved} & Abs.\ gap \\
     \midrule
     50 & 5 & \textbf{100\%} & \textbf{2.12e-08} & \textbf{100\%} & 3.53e-05 & \textbf{100\%} & 1.76e-07 & \textbf{100\%} & 2.56e-07 & \textbf{100\%} & 2.93e-07 & \textbf{100\%} & 3.15e-07 \\
     50 & 12 & \textbf{100\%} & 4.19e-07 & \textbf{100\%} & \textbf{1.07e-07} & \textbf{100\%} & 3.22e-07 & \textbf{100\%} & 2.45e-07 & \textbf{100\%} & 2.30e-07 & \textbf{100\%} & 2.67e-07 \\
     60 & 6 & \textbf{100\%} & \textbf{1.17e-08} & \textbf{100\%} & 5.16e-08 & \textbf{100\%} & 1.77e-07 & \textbf{100\%} & 3.02e-07 & \textbf{100\%} & 1.80e-07 & \textbf{100\%} & 3.72e-07 \\
     60 & 15 & \textbf{100\%} & 7.29e-07 & \textbf{100\%} & 7.87e-07 & \textbf{100\%} & \textbf{4.94e-07} & 80\% & 5.74e-05 & 80\% & 7.66e-05 & 80\% & 1.93e-04 \\
     80 & 8 & \textbf{100\%} & 4.33e-07 & \textbf{100\%} & 1.94e-07 & \textbf{100\%} & 2.47e-07 & \textbf{100\%} & 6.15e-07 & \textbf{100\%} & 4.38e-07 & \textbf{100\%} & \textbf{8.11e-08} \\
     80 & 20 & 0\% & 6.43e-03 & 0\% & \textbf{5.79e-03} & 0\% & 6.04e-03 & 0\% & 6.48e-03 & 0\% & 6.84e-03 & 0\% & 7.29e-03 \\
     100 & 10 & \textbf{100\%} & 5.78e-07 & \textbf{100\%} & 4.78e-07 & \textbf{100\%} & 5.50e-07 & \textbf{100\%} & \textbf{4.12e-07} & \textbf{100\%} & 7.10e-07 & \textbf{100\%} & 7.40e-07 \\
     100 & 25 & 0\% & 1.06e-02 & 0\% & \textbf{9.85e-03} & 0\% & 9.87e-03 & 0\% & 1.10e-02 & 0\% & 1.13e-02 & 0\% & 1.16e-02 \\
     120 & 12 & 60\% & 1.48e-04 & \textbf{100\%} & \textbf{9.61e-07} & 80\% & 1.66e-05 & 80\% & 8.47e-05 & 80\% & 8.64e-05 & 40\% & 1.93e-04 \\
     120 & 30 & 0\% & 1.38e-02 & 0\% & \textbf{1.26e-02} & 0\% & 1.31e-02 & 0\% & 1.38e-02 & 0\% & 1.39e-02 & 0\% & 1.47e-02 \\
     \bottomrule
    \end{tabular}%
    }
    \caption{\protect\revII{Abs.\ gap is the mean absolute gap for different Frank-Wolfe tolerances on D-optimal design with correlated data.
    The best result in each row is shown in bold.}}
    \label{tab:final_abs_gap_d_correlated}
\end{table}

\begin{table}[htbp]
    \centering
    \resizebox{\linewidth}{!}{%
    \begin{tabular}{ccrrrrrrrrrrrr}
    \toprule
    \# exp. & \# par. & \multicolumn{2}{c}{5e-2} & \multicolumn{2}{c}{1e-2} & \multicolumn{2}{c}{5e-3} & \multicolumn{2}{c}{1e-3} & \multicolumn{2}{c}{5e-4} & \multicolumn{2}{c}{1e-4} \\
    \cmidrule(lr){3-4}\cmidrule(lr){5-6}\cmidrule(lr){7-8}\cmidrule(lr){9-10}\cmidrule(lr){11-12}\cmidrule(lr){13-14}
     &  & \makecell{\%\\solved} & Abs.\ gap & \makecell{\%\\solved} & Abs.\ gap & \makecell{\%\\solved} & Abs.\ gap & \makecell{\%\\solved} & Abs.\ gap & \makecell{\%\\solved} & Abs.\ gap & \makecell{\%\\solved} & Abs.\ gap \\
     \midrule
     50 & 5 & \textbf{100\%} & \textbf{0.00e+00} & \textbf{100\%} & \textbf{0.00e+00} & \textbf{100\%} & \textbf{0.00e+00} & \textbf{100\%} & \textbf{0.00e+00} & \textbf{100\%} & \textbf{0.00e+00} & \textbf{100\%} & 1.35e-08 \\
     50 & 12 & \textbf{100\%} & 6.09e-07 & \textbf{100\%} & 5.90e-07 & \textbf{100\%} & 5.01e-07 & \textbf{100\%} & \textbf{4.16e-07} & \textbf{100\%} & 6.15e-07 & \textbf{100\%} & 5.28e-07 \\
     60 & 6 & \textbf{100\%} & 1.92e-07 & \textbf{100\%} & 4.84e-09 & \textbf{100\%} & \textbf{0.00e+00} & \textbf{100\%} & 3.51e-07 & \textbf{100\%} & 1.93e-07 & \textbf{100\%} & 2.12e-07 \\
     60 & 15 & \textbf{40\%} & 2.22e-03 & \textbf{40\%} & \textbf{1.92e-03} & 20\% & 2.15e-03 & 20\% & 2.52e-03 & 20\% & 2.64e-03 & 0\% & 2.97e-03 \\
     80 & 8 & \textbf{100\%} & 4.16e-07 & \textbf{100\%} & \textbf{1.41e-07} & \textbf{100\%} & 2.94e-07 & \textbf{100\%} & 1.67e-07 & \textbf{100\%} & 1.62e-07 & \textbf{100\%} & 3.64e-07 \\
     80 & 20 & 0\% & 7.64e-03 & 0\% & \textbf{7.17e-03} & 0\% & 7.33e-03 & 0\% & 7.80e-03 & 0\% & 8.08e-03 & 0\% & 8.30e-03 \\
     100 & 10 & \textbf{100\%} & 9.71e-07 & \textbf{100\%} & 9.33e-07 & \textbf{100\%} & 9.36e-07 & \textbf{100\%} & 9.05e-07 & \textbf{100\%} & \textbf{8.86e-07} & 80\% & 3.06e-05 \\
     100 & 25 & 0\% & 1.44e-02 & 0\% & \textbf{1.35e-02} & 0\% & 1.38e-02 & 0\% & 1.46e-02 & 0\% & 1.47e-02 & 0\% & 1.50e-02 \\
     120 & 12 & 0\% & 1.03e-03 & \textbf{40\%} & 4.74e-04 & \textbf{40\%} & \textbf{4.29e-04} & \textbf{40\%} & 5.11e-04 & 20\% & 5.75e-04 & 20\% & 7.43e-04 \\
     120 & 30 & 0\% & 1.90e-02 & 0\% & 1.76e-02 & 0\% & \textbf{1.74e-02} & 0\% & 1.84e-02 & 0\% & 1.90e-02 & 0\% & 2.02e-02 \\
     \bottomrule
    \end{tabular}%
    }
    \caption{\protect\revII{Abs.\ gap is the mean absolute gap for different Frank-Wolfe tolerances on D-optimal design with independent data.
    The best result in each row is shown in bold.}}
    \label{tab:final_abs_gap_d_independent}
\end{table}

\end{document}